\documentclass[12pt]{amsart}

\def\@typesizes{%
       \or{5}{6.5}\or{6}{7.5}\or{7}{8.5}\or{8}{11}\or{9}{12}%
       \or{10}{13}
       \or{\@xipt}{14}\or{\@xiipt}{15}\or{\@xivpt}{18}%
       \or{\@xviipt}{20}\or{\@xxpt}{24}}

\usepackage{amsthm}
\usepackage{amsmath}
\usepackage{amssymb}
\usepackage{graphicx}
\usepackage{enumerate}
\allowdisplaybreaks
\numberwithin{equation}{section}
\numberwithin{figure}{section}

\theoremstyle{plain}
\newtheorem{theorem}{ Theorem}[section]
\newtheorem{proposition}[theorem]{ Proposition}
\newtheorem{lemma}[theorem]{ Lemma}
\newtheorem{corollary}[theorem]{ Corollary}
\newtheorem{example}[theorem]{ Example}
\newtheorem{remark}[theorem]{ Remark}
\newtheorem{definition}[theorem]{ Definition}
\newtheorem{conjecture}{ Conjecture}

\def\BET{\begin{theorem}}
\def\ENT{\end{theorem}}
\def\BEP{\begin{proposition}}
\def\ENP{\end{proposition}}
\def\BEL{\begin{lemma}}
\def\ENL{\end{lemma}}
\def\BEC{\begin{corollary}}
\def\ENC{\end{corollary}}
\def\BEE{\begin{example} \rm}
\def\ENE{\end{example}}
\def\BER{\begin{remark} \rm}
\def\ENR{\end{remark}}
\def\BED{\begin{definition} \rm}
\def\END{\end{definition}}
\def\BECJ{\begin{conjecture}}
\def\ENCJ{\end{conjecture}}

\def\bea{\begin{eqnarray}}
\def\eea{\end{eqnarray}}

\def\beq{\begin{equation}}
\def\eeq{\end{equation}}

\def\beal{\begin{align*}}

\def\eeal{ \end{align*} }

\def\roweq{\nonumber \\ &=& }
\def\rowleq{\nonumber \\  & \leq & }

\def\rowpl{\nonumber \\ & \ \ + & }
\def\rowmi{\nonumber \\ & \ \ - & }

\def\bfB{{\bf B}}

\def\bfE{{\bf E}}

\def\bfM{{\bf M}}

\def\bfW{{\bf W}}

\def\bbC{{\mathbb C}}

\def\bbN{{\mathbb N}}

\def\bbR{{\mathbb R}}

\def\cA{{\mathcal A}}
\def\cB{{\mathcal B}}

\def\cF{{\mathcal F}}

\def\cH{{\mathcal H}}

\def\cT{{\mathcal T}}
\def\cU{{\mathcal U}}
\def\cV{{\mathcal V}}
\def\cW{{\mathcal W}}
\def\cX{{\mathcal X}}
\def\cY{{\mathcal Y}}

\newcommand{\diag}{\mathop{\rm diag}\nolimits}

\begin{document}

\title{Bands in the spectrum of a periodic elastic waveguide }

\author{F.L. Bakharev}
\author{J. Taskinen}

\thanks{The first named author was supported by the St. Petersburg State University
grant  6.42.1372.2015, by Russian Foundation for Basic Research (Grant 15–01–02175), 
by the Chebyshev Laboratory - RF Government grant 11.G34.31.0026, and by JSC "Gazprom Neft". 
The second named authors was  supported by grants from
the Magnus Ehrnroot Foundation and the V\"ais\"al\"a Foundation of the
Finnish Academy of Sciences and Letters.} 

\address{
Chebyshev Laboratory, St. Petersburg State University, 14th Line, 29b, Saint Petersburg, 199178 Russia}

\address{University of Helsinki, Department of 
Mathematics and Statistics,
P.O. Box 68, FI-00014 Helsinki, Finland. }

\maketitle

{\small

{\bf Abstract}:
We study the spectral linear elasticity problem in an unbounded 
periodic waveguide, which consists of a sequence of identical bounded
cells connected by thin ligaments of diameter of order $ h >0$. 
The essential spectrum of the problem is known to have band-gap structure.
We derive asymptotic formulas for the position of the spectral bands and
gaps, as $h \to 0$. 
}

\bigskip

\section{Introduction}
\label{sec1}
\begin{figure}
\begin{center}
\includegraphics[width=10cm]{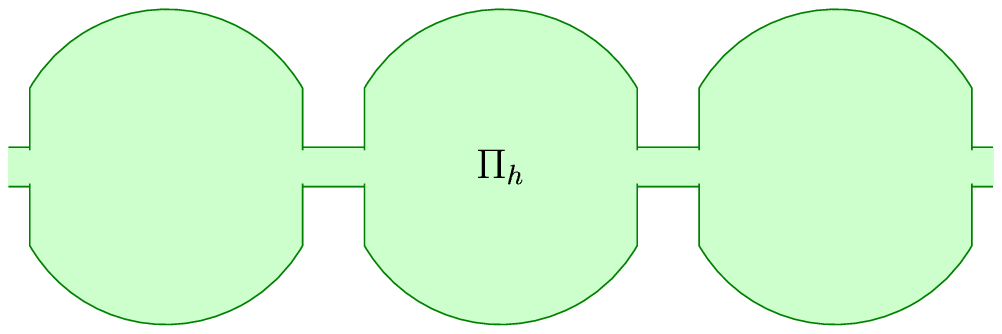}
\end{center}
\caption{The waveguide $\Pi_h$}
\label{fig2}
\end{figure}

We study the essential spectrum 
of the  linearized elasticity system with 
tract\-ion--free boundary conditions in unbounded periodic waveguides
denoted by $\Pi_h$. 
The waveguide, see Fig.\,\ref{fig2},
consists of infinitely many identical, translated bounded cells
connected with small cylindrical ligaments, the length and radius of cross-section 
of which are both proportional to a small parameter $h>0$ so that the volume of 
the ligament is $O(h^3)$. 
A number of papers (e.g \cite{BaRuTa}, \cite{naruta2}, \cite{nata})
has been devoted to geometrically similar waveguides consisting
of arrays of macroscopic cells connected with thin structures, and,  using 
rigorous perturbation arguments, 
the existence of gaps in the essential spectra has been detected. This has been
done for  elliptic boundary problems in elasticity, linear water-wave theory, 
piezo-electricity etc. 

In this geometric setting the emerging of gaps is explained by that for small $h$, the 
problem can be seen as a perturbation of a "limit spectral problem" ($h=0$) on 
a bounded domain, which consists of a single cell $\varpi_0$. The spectrum of 
such a problem is in general a sequence $(\lambda_k)_ {k=1}^\infty$ 
of eigenvalues , and the spectral bands of the original
problem are situated "close" to the eigenvalues $\lambda_k$. To analyse 
this closeness and its dependence on $h$ becomes a mathematical challenge: 
if that can be done accurately enough, one finds that disjoint eigenvalues
correspond to spectral bands with a gap in  between. In fact, this scheme
can work only for a finite number of the lowest eigenvalues in the sense
that for each fixed small $h$, at most a finite number of gaps can be
found.

The purpose of this work is to refine the existing results by proving more
accurate estimates than before for the end-points of spectral bands and gaps. 
For example, in  \cite{naruta2} it was shown that for small $h$ and also
$k$, the $k$th spectral band is situated within a distance $Ch$ from 
$\lambda_k$, the $k$th eigenvalue of the limit problem, where $C>0$ is
a constant depending on the shape of the cells and on the physical constants
of the elastic material, but not on the size of the
small ligaments. In comparison, we shall find here an asymptotic formula for the
spectral bands $\Upsilon_k$. Namely, by the Floquet-Bloch-theory of periodic problems, the 
bands are formed by eigenvalues $\Lambda_k^h(\eta)$ of the "model problem" 
depending on the  parameter $\eta\in [0,2 \pi)$. It is in fact well known that the
essential spectrum of the original problem (later \eqref{prob-1}--\eqref{prob-2})
equals 
\begin{equation}
\sigma_{\rm ess} = \bigcup_{k=1}^\infty \Upsilon_k^h \ , \ \ 
\Upsilon_k^h = \{  \Lambda_k^h (\eta) \, : \,
\eta \in [0,2 \pi) \} .
\label{eq1.1}
\end{equation}
In Theorem \ref{th6.2}, see \eqref{6.12},
\eqref{1.98}, we determine the 
first order (in $h$) correction term for the difference of $\Lambda_k^h(\eta)$
and $\lambda_k$. The result includes the following claim: 

{\it For all $k$ we have the estimate
\bea
\big| \Lambda_k^h(\eta) -  ( \lambda_k + h\Lambda_k'(\eta)) 
\big| \leq C_k h^{3/2} \ \ \forall h > 0\ ,
\ \eta \in [0,2\pi) , \label{1.99}
\eea
where  
\bea
\Lambda_k'(\eta) = 2 a \lambda_k  + \big(A_k + e^{i\eta} B_k\big)^\top
\bfM^+ \big(A_k + e^{-i\eta} B_k\big) , \label{1.98p}
\eea
the column vectors $A_k$, $B_k \in \bbR^3$ and the positive definite matrix $\bfM^+ \in \bbR^{3 \times3}$ do not depend on $h$ or $\eta$, and the number $a$ equals $0$ or $1$ (according to  the domain, see Section \ref{sec2.1}). }

The quantities $A_k$, $B_k$, $\bfM^+$ will be determined in Sections \ref{sec5.2}, \ref{sec5.3}. Anyway, the coefficient $\Lambda_k'(\eta)$ of the correction
term depends only on the geometry of the limit problem.   
(It might thus be desirable to provide numerical experiments
on the effect of this coefficient and the geometry of the limit domain to the
appearance of spectral gaps, but we do not provide such information here.)
However, information on the position and length of $\Upsilon_k$  obtained here is
much more precise than before, since until now there has not even existed
a criterion to distinguish, if the band has positive length or it consists of only a single point, which is an eigenvalue of infinite multiplicity. In many cases 
the vector $B_k$ is nonzero, and it follows from \eqref{1.98p} that length of the band
$\Upsilon_k^h$ is positive.

In addition, we find an asymptotic representation
for the corresponding eigenfunctions  $U_k^h$  of the model problem, which includes the 
leading terms both near the junctions of the periodicity cells and at a
distance of these points. Since it  is difficult to describe this 
without a number of definitions, we refer to Theorem \ref{th6.2} for details.

In literature spectral gaps in essential spectra for scalar equations and Maxwell's system 
in infinite periodic media have been considered  in many papers,  
see for example \cite{Kuchreview},  \cite{Filon}, \cite{Fried}, \cite{Green}, 
\cite{HeLi}, \cite{Zh}. For results on gaps in (quasiperiodic, unbounded) waveguides we mention the  papers  \cite{BaNaRu}, \cite{CaNaPe}, \cite{FriedSol}, \cite{NaMatZam};
for an  approach  based on  parameter--dependent Korn--type inequalities, see
\cite{CaNaZAMM}, \cite{naCRAS}, \cite{naDAN}, \cite{na2009}, \cite{naJVMMF}. 
A comparison of the present work with  the paper \cite{naruta2} was already
presented above. We just mention that  the method of \cite{naruta2} is based on the max-min 
principle for eigenvalues. 
We finally also mention the paper \cite{BaRuTa}, which contains an analysis of spectral bands much like in the present work, but in the more simple 
setting of the linear water wave equation.

As for the structure of this paper, we present in Section \ref{sec2} the geometry
of the waveguide $\Pi_h$, the  formulation of the linear spectral elasticity problem and  its 
variational formulation.  The parameter dependent problem arising from the FBG-transform
is presented in Section \ref{sec2.2}, together with the  variational formulation in a special Sobolev-type Hilbert space
$\cH^{h,\eta}$. The limit problem is studied in Section \ref{sec2.3}, and
additional technical devices are introduced in Section \ref{sec2.4}. 
Section \ref{sec5} contains the formal asymptotic analysis, in particular the construction 
of the asymptotic ans\"atze for the eigenfunctions of the model problem.
Section \ref{sec6.1} contains the main result, and proof is given in Sections \ref{sec6.2}--\ref{sec6.3}. One point of the proof needs the existence of spectral 
gaps, which can be obtained by an adaptation of the method \cite{naruta2}; this
is presented in the appendix, Section \ref{secNAR1}.

Acknowledgement. The authors want to thank Prof. Sergei A. Nazarov for many discussions on the topic of this work.

\section{Problem formulation, band-gap spectrum and limit problem}
\label{sec2}

\subsection{Spectral elasticity problem in the waveguide.}
\label{sec2.1}
The problem domain  is a periodic waveguide $\Pi_h$, Fig.\,\ref{fig2}, depending on a small geometric
parameter $h \in (0,\frac1{10}]$ and consisting of infinitely many disjoint translated copies of a bounded cell  
connected by  thin cylinder.  Let us describe the exact definition. We denote by $\varpi_\bullet$  a bounded domain in ${\mathbb R}^3$ containing the points $P^\pm=(0,0,\pm 1/2)$ in its interior, such that 
\begin{equation}
\label{varpi}
\varpi_0=\{x=(x_1,x_2,x_3)=(y,z)\in \varpi_\bullet \,:\, y\in {\mathbb R^2}, |z|<1/2\}
\end{equation}
is a domain with Lipschitz boundary $\partial \varpi_0$. Notice that in a neighbourhood of the 
points $P^\pm$ the boundary  $\partial  \varpi_0$ is "flat": it consists of an open subset of the plane
$\{ (y_1,y_2, \pm 1/2)\}$.
  
We also denote by  $\theta\subset {\mathbb R}^2$ a bounded domain which has Lipschitz boundary $\partial \theta$ and which is star-shaped with respect to the point $(0,0) \in \theta$. (Any point of 
$\theta$ can be joined with $(0,0)$ by a line segment running in $\theta$.)
We also denote  for $h>0$
\begin{equation*}
\theta_h=h\theta=\{ h y\in{\mathbb R}^2:y\in \theta\}.
\end{equation*}
By a thin cylinder  we mean the set $\Theta_h=\theta_h\times {\mathbb R}$.

\begin{figure}
\begin{center}
 \includegraphics[width=10cm]{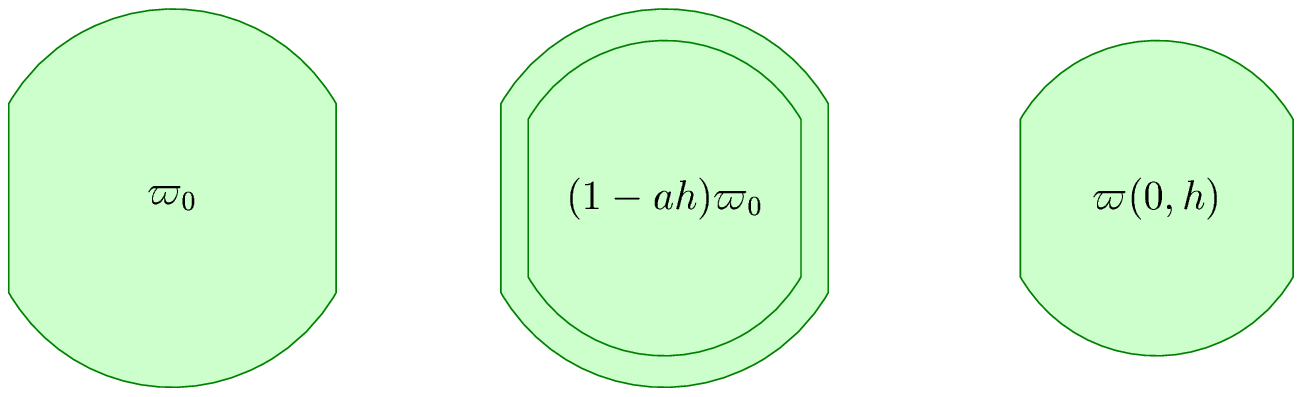}
\end{center}
\caption{Transformation of the set $\varpi_0$}
\label{fig1}
\end{figure}

\begin{figure}
\begin{center}
\includegraphics[height=4cm,width=10cm] {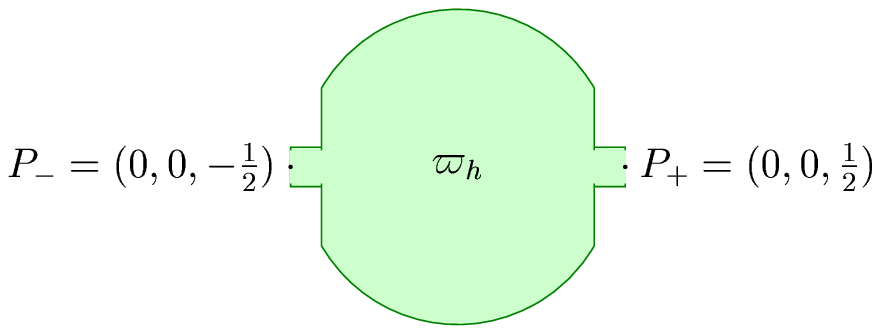} 
\end{center}
\caption{Periodicity cell $\varpi_h$}
\label{fig15}
\end{figure}

We fix a number $a$ to be either 0 or 1 (see below for a discussion) 
and denote by  $a_h$ the number $( 1-ah)^{-1} \in \bbR$ and also the multiplier
$
a_h : x \mapsto (1-ah)^{-1} x $ 
in any space $\bbR^n\ni x$, $n \in \bbN$. The meaning will be clear from the context. 
Let us  define a scaled, translated  family of bodies 
\begin{equation}
\label{varpi(j,h)}
\varpi(j,h)=\{x : a_h(y,z-j)\in \varpi_0 \} \ \ \mbox{with} 
\  \varpi(h) := \varpi(0,h) = a_h^{-1} (\varpi_0),
\end{equation}
see Fig.\,\ref{fig1}.  The periodicity cell  $\varpi_h \supset \varpi(h)$, 
Fig.\,\ref{fig15}, is the set
\begin{equation}
\label{varpi_h}
\varpi_h=\varpi(h)\cup \Big( \theta_h\times \big(-\frac{1}{2},\frac{1}{2}\big)\Big),
\end{equation}
in other words, $\varpi_h$ is the disjoint union of $\varpi(h)$ and the small
sets
\begin{equation}
\label{Theta_h}
\Theta_h^-=\theta_h\times (-1/2,-1/2+ah/2],\quad \Theta_h^+=\theta_h\times (1/2-ah/2,1/2].
\end{equation} 
Finally, the  waveguide $\Pi_h$, Fig.\,\ref{fig2}, 
is  defined as
\begin{equation}
\label{Pi_h}
\Pi_h=\Theta_h\cup \bigcup_{j\in {\mathbb Z}} \varpi(j,h) .
\end{equation}
Note that when $h=0$, the waveguide $\Pi_h$ becomes  a union of disconnected 
sets.

If $a=0$ above, then of course $a_h =1$  is independent of $h$, and the cells
$\varpi(j,h)$ contact each other: the cylinder $\Theta_h$  loses its role and the cells are 
connected by holes, "apertures" of diameter $h$  in their boundaries. In this case
also the volume of the  periodicity cells is independent of $h$. 
The subsequent calculations are concentrated on treating the case $a=1$. In many
cases the calculations are unnecessarily complicated for $a=0$, but that case  
can be considered as a by-product; some details will be omitted, and the reader is 
asked to keep it also in mind.  If $a=1$, the diameter of the connecting
ligaments is of order $h$ and the volume is proportional $h^3$. This means a certain
difference to the related  paper \cite{naruta2}.

We consider the spectral elasticity problem  in $\Pi_h$ 
written in matrix form, cf. \cite{Lekh}, \cite{Nabook}. This is formulated for an unknown 
$\bbR^3$-valued vector function
$u(x) $,  
which describes the displacement vector of the given material. 
The system contains a first order differential operator
matrix $D(\nabla_x)$, where $\nabla_x$ denotes the gradient with respect to the
variable $x$ and 
\begin{equation}
\nonumber
D(x)=
\left(
\begin{array}{cccccc}
x_1&0&2^{-1/2}x_2&2^{-1/2}x_3&0&0\\
0&x_2&2^{-1/2}x_1&0&2^{-1/2}x_3&0\\
0&0&0&2^{-1/2}x_1&2^{-1/2}x_2&x_3
\end{array}
\right)^{\top}\,, \quad x=(x_1,x_2,x_3)^\top\,,
\end{equation}
and a matrix $A$ of dimension $6 \times 6$, describing elastic moduli. The matrix $A$  is assumed to be positive definite and, for technical simplicity, constant.
The problem  can be written as 
\begin{eqnarray}
\label{prob-1}                                                            
L_x 
u(x) & = & 
\lambda u(x)\ , \ \  x\in \Pi_h,\\
\label{prob-2}
{N_x^h} 
u(x)& = & 0 \  \mbox{for a.e.} \  x\in \partial \Pi_h\, ,
\end{eqnarray}
where
\bea
L_x  = D(-\nabla_x)^\top A D(\nabla_x) \ , \ \
{N_x^h} = D(\nu(x))^\top AD(\nabla_x) \ , 
\label{4.77u}
\eea 
$\nu$ is the outward normal vector defined for almost all points of the
Lipshitz surface $\partial \Pi_h$ and $\lambda$ is a spectral parameter.
The boundary operator ${N_x^h}$ depends on $h$ via the  domain
$\Pi_h$.  (Later  on, $\nu$ will always denote the outward normal
of a boundary, which will be clear from the context.) 
Also,  $L_x$ is just a constant coefficient partial 
differential operator containing only second order terms.
In general, $\lambda$ could be multiplied in \eqref{prob-1} by a fixed function describing the 
material density, but we assume this function to be equal to 1 for simplicity. 

The variational formulation of the spectral problem \eqref{prob-1}-\eqref{prob-2} reads as
\begin{equation}
\label{variat}
a (u,  v; \Pi_h )=\lambda(u, v)_{\Pi_h}, \quad v\in H^1(\Pi_h)^3.
\end{equation}
Here, given a regular enough domain $\Omega \subset \bbR^3$, we use the notation
\bea
& & a(f,g; \Omega) = (A D(\nabla_x) f , D(\nabla_x) g)_{\Omega}
\ , \ \mbox{and} \ \  (f,g)_{\Omega}=\sum_{j=1}^k \int\limits_{\Omega}f \overline{g} dx
\label{1.31}
\eea
for  the usual (complex valued) inner product in $L^2(\Omega)^k$, $k =1,2,3$;
the latter notation will also be used for domains in $\bbR^2$.
We denote the standard Sobolev space of first order on $\Omega$ by $H^1(\Omega)$.
The bilinear form \eqref{1.31}
is positive and closed in the Sobolev space $H^1(\Pi_h)^3$ and consequently (see \cite{BiSo})
our problem can be rewritten as an abstract operator equation
$
\cT^h u=\lambda u
$,
where $\cT^h$ is an unbounded, self-adjoint densely defined operator in the Hilbert space $L^2(\Pi_h)^3$ and thus
the spectrum $\sigma(\cT^h) \ni \lambda$ is a subset of $\overline{{\mathbb R}_+}=[0,+\infty)$.
The embedding $H^1(\Pi_h)\subset L^2(\Pi_h)$ is not compact due to the  unboundedness of the domain $\Pi_h$,
hence, the essential spectrum $\sigma_{\rm{ess}}(\cT^h)$ is not empty (see \cite{BiSo}, Th.10.15).
Finally, spectral concepts of the problem \eqref{prob-1}--\eqref{prob-2} are 
defined with the help
of the operator $\cT^h$. In particular, $\sigma_{\rm{ess}} := \sigma_{\rm{ess}}(\cT^h)$.

\subsection{Floquet-Bloch-Gelfand transform and essential spectrum.}
\label{sec2.2}
To analyse the band-gap structure spectrum  of  the problem 
\eqref{prob-1}--\eqref{prob-2} we use the FBG-transform
\begin{equation}
\nonumber 
v(y,z)\mapsto\ V(y,z ; \eta)=\frac{1}{\sqrt{2\pi}}\sum_{j\in \mathbb{Z}}
\exp(-i\eta j)v(y,z+j),
\end{equation}
where $(y,z)\in \Pi_h$ on the left, while $\eta\in[0,2\pi)$ and
$(y,z)\in\varpi_h$ on the right. As well known, this operator establishes an
isometric isomorphism between Lebesgue spaces
$
L^2(\Pi_h)^3$ and $L^2(0,2\pi; L^2(\varpi_h)^3),
$
where $L^2(0,2\pi;B)$ is the Lebesgue space of functions with values in the
Banach space $B$, endowed with the norm
$
\|V;L^2(0,2\pi; B)\|^2= \int_0^{2\pi}\|V(\eta);B\|^2d\eta \,.
$
The FBG-transform is also an isomorphism from the Sobolev space 
$H^1(\Pi_h)^3$ onto 
$L^2(0,2\pi; H^1_\eta(\varpi_h)^3)$ and  from $H^2(\Pi_h)^3$ onto 
$L^2(0,2\pi; H^2_\eta(\varpi_h)^3)$. Here, for a fixed $h$ and $\eta$, the  space $H^2_{\eta}(\varpi_h)$ 
is the space of Sobolev functions $f$ on $\varpi_h$,
which satisfy quasiperiodicity conditions 
\begin{eqnarray}
\label{qua-1} \ \ 
f(y,1/2) & = & e^{i\eta} f (y,-1/2) \ , \ \  
\partial_{z} f (y,1/2) =  e^{i\eta}\partial_{z} f (y,-1/2)\ , \ \  y\in \theta_h.
\end{eqnarray}
Similarly, the  space $H^1_{\eta}(\varpi_h)$ consists of $H^1$-functions
satisfying the condition \eqref{qua-1} only. 
In the following we denote by $\cH^{h,\eta}$ the space $H^1_{\eta}(\varpi_h)^3$
endowed with the norm $\Vert f ; \cH^{h,\eta} \Vert$
coming from the inner product 
\bea
\langle f,g \rangle_{\cH^{h,\eta}} = 
(f,g)_{\varpi_h} 
+ \big(A D(\nabla_x) f, D(\nabla_x) g \big)_{\varpi_h} ; \label{2.12a}
\eea
notice that  the last form on the right is Hermitian and 
positive on $H^1_\eta(\varpi_h)^3$.
Obviously, $\Vert f ; \cH^{h,\eta} \Vert \leq C \Vert f ; H^1 (\varpi_h) \Vert$ for some constant independent of $h$ and $\eta$, since the matrix $A$ was assumed to be constant. 

\BEL
\label{lem2.1}
The norm $\Vert \cdot ; \cH^{h,\eta} \Vert $ is equivalent to $ \Vert \cdot ; H^1 (\varpi_h) \Vert$
with constants independent of $h > 0$ and $\eta \in [0, 2 \pi)$. 
\ENL

Proof. It is enough to show that the Korn inequality
\bea
 \Vert f ; H^1 (\varpi_h) \Vert^2 \leq C\big(  \Vert f ;  L^2 (\varpi_h) \Vert^2 
 + \big(A D(\nabla_x) f, D(\nabla_x) f \big)_{\varpi_h} \label{2.12cc}
\eea
holds  for Sobolev-functions $f$ in such a way that the constant $C>0$ can be chosen 
independently of $h$, $\eta$. It is quite obvious that replacing $\varpi_h$ by 
$\varpi(h)$ (see \eqref{varpi(j,h)})  in \eqref{2.12cc}, the corresponding
inequality would hold with constant independent of $h$, $\eta$. But the same
is true also, when $\varpi_h$ is replaced by the 
two sets  $\Theta_h^\pm$, see \eqref{Theta_h}.  
This follows from the result \cite{KoOl89},
Th.\,1, Th.\,2. For example the set $\Theta_h^+ $
has diameter $c_1h$ and it is  in the terminology of the citation 
starshaped with respect to a ball with center in  $(0,0,1/2 - a h/4 )$ and radius $c_2 h$, where 
$c_1$,  $c_3$ do not depend on $h$ (see the assumptions on 
$\theta$ in Section \ref{sec2.1}). Then, \eqref{2.12cc} follows by combining these facts.
 \ \ $\Box$

\bigskip

Using the FBG-transform the problem \eqref{prob-1}-\eqref{prob-2} turns into a model
problem on the periodicity cell for the unknown $U = U(x;\eta)$,
\begin{eqnarray}
\label{prob-2-1}                                                            
& & \overline{D(-\nabla_y, -\partial_z-i\eta)}^\top A 
D(\nabla_y, \partial_z+i\eta) U  =  \Lambda  U \ \ \  \mbox{in} \  \varpi_h,\\
\label{prob-2-2}
& & \overline{D(\nu(x))}^\top A D(\nabla_y, \partial_z+i\eta) U  =    0 \ \ \  \mbox{ a.e. in} \ 
\partial \varpi_h\setminus (\theta_h^-\cup\theta_h^+)\,, \\
\label{prob-2-3}
& & U(y,1/2;\eta) =  e^{i\eta}U(y,-1/2;\eta), \quad y\in \theta_h\,,\\
\label{prob-2-4}
& & \partial_{z}U(y,1/2;\eta) =  e^{i\eta}\partial_{z}U(y,-1/2;\eta), \quad y\in \theta_h\,,
\end{eqnarray}
where the overline denotes complex conjugation and $\Lambda$ is a spectral parameter. In the weak form this amounts  to finding $0 \not= U \in 
\cH^{h,\eta}$ and $\Lambda \in \bbC$ with
\begin{equation}
\label{variat-2}
a (U,V ; \varpi_h )=\Lambda(U,V)_{\varpi_h}
\end{equation}
for all $V\in 
\cH^{h,\eta}$.  The problem 
\eqref{prob-2-1}-\eqref{prob-2-4} can be associated with a self-adjoint, 
positive and compact operator  
$
\mathcal{B}^{h,\eta}: \cH^{h,\eta} \to  \cH^{h,\eta}.
$
We define $\cB^{h,\eta}$ in a standard way by requiring that the identity
\begin{equation}
\label{3.10}
\langle \cB^{h,\eta} U,V  \rangle_{\cH^{h,\eta}}=(U,V)_{\varpi_h}
\end{equation}  
holds for all $U,V \in \cH^{h,\eta}$.
The problem \eqref{variat-2} is then  equivalent to the spectral problem
$
\cB^{h,\eta}U =\mu U
$
with a new spectral parameter $\mu=(1+\Lambda)^{-1}$. The spectrum
of $\cB^{h,\eta}$ consists of a 
decreasing sequence $\{\mu_k^h(\eta)\}_{k\geq 1}$ of eigenvalues and the point 0 of the essential  spectrum. As a consequence, the spectrum of the problem 
\eqref{prob-2-1}-\eqref{prob-2-4}
can be presented as the eigenvalue sequence (counting multiplicities)
\begin{equation}
\label{Lambda}
0\leq \Lambda_1^h(\eta)\leq \Lambda_2^h(\eta)\leq \ldots \leq \Lambda^h_p(\eta)\leq\ldots\to+\infty.
\end{equation}
We denote the corresponding eigenfunctions by $U^h_p  \in H^1_\eta(\varpi_h)^3$,
where  the dependence on $\eta$ will usually not be shown. We require the orthogonality property
\begin{equation}
\label{ortho-U}
(U^h_p,U^h_q)_{\varpi_h}=\delta_{p,q}, \quad p,q=1, 2, \ldots . 
\end{equation}
By \cite{Gel}, \cite{Kuch}, \cite{na17}, \cite{NaSpec}, \cite{NaPl}, 
Theorem 3.4.6,  \cite{naruta2}, Theorem 2.1), for example, 
a number $\Lambda$ belongs to the resolvent set or the  discrete spectrum of  $\cT^h$ (end of 
Section \ref{sec2.1}),  if and only if it does not coincide with 
$\Lambda_p^h(\eta)$ for any $\eta \in [0,2\pi]$ and $p$.  Hence,
the essential spectrum of $\cT^h$
and thus also of the original problem 
\eqref{prob-1}--\eqref{prob-2} have band-gap structure \eqref{eq1.1},
\begin{equation}
\label{band-gap}
\sigma_{{\rm ess}}(\cT^h)=\bigcup_{p=1}^\infty \Upsilon_p^h,
= \bigcup_{p=1}^\infty
\{\Lambda_p^h(\eta) : \eta\in[0,2\pi)\}; 
\end{equation}
where the spectral bands $\Upsilon_p^h$ are closed intervals (possibly single points).

\subsection{Spectrum of the limit model problem}
\label{sec2.3}
When $h \to 0$, the cylinder $\Theta_h$ turns into a negligible set 
and the quasi-periodicity conditions \eqref{prob-2-3}-\eqref{prob-2-4} 
lose their meaning. The model problem \eqref{prob-2-1}-\eqref{prob-2-4} 
turns into  the so called limit problem 
on the isolated elastic body $\varpi_0$,
\begin{eqnarray}
\label{prob-3-1}
&&D(-\nabla_x)^\top A D(\nabla_x) u=\Lambda^0 u \ \ \mbox{in} \ \varpi_0,\\
\label{prob-3-2}
&&D(\nu(x))^\top A D(\nabla_x) u= 0 \ \  \mbox{ a.e.\,in} \  \partial \varpi_0 \,,
\end{eqnarray}
where $u$ is the unknown function on $\varpi_0$ and $\Lambda^0$ is a spectral
parameter. 
The problem has for any  $\eta \in (0,2\pi]$  the same eigenvalues as the case $\eta=0$, so we can restrict to this single limit case, see \cite{naruta2} for some more details.
We again denote by  $a(u,v; \varpi_0)$ the sesquilinear form corresponding to the problem \eqref{prob-3-1}-\eqref{prob-3-2}. 
It is 
positive 
and closed on $H^1(\varpi_0)^3$. We also define the space $\cH^0$ (with no
quasiperiodicity conditions) and the self-adjoint, positive, compact
operator $\cB^0 : \cH^0 \to \cH^0$, associated to the problem
\eqref{prob-3-1}--\eqref{prob-3-2}, analogously to $\cH^{h,\eta}$ and  
$\cB^{h,\eta}$ of the previous section. 
Since $\varpi_0$ is a bounded Lipschitz domain, the limit problem has an eigenvalue sequence
\begin{equation}
\label{lambda}
0=\lambda_1=\ldots=\lambda_6<\lambda_7\leq\lambda_8\leq\ldots\leq\lambda_p\leq\ldots \to+\infty.
\end{equation}
Here the $\lambda_1$, \ldots, $\lambda_6$ are the eigenvalues of the rigid motions
including three translations and rotations. For every $p$, let
$u_p$ be the eigenfunction corresponding to $\lambda_p$, normalized in
$L^2(\varpi_0)^3$ by 
\begin{equation}
\label{ortho-u}
(u_p,u_q)_{\varpi_0}=\delta_{p,q}, \quad p,q = 1, 2,\ldots
\end{equation}

We also state a scaled version of the problem \eqref{prob-3-1}--\eqref{prob-3-2}
(see \eqref{varpi(j,h)}):
\begin{eqnarray}
\label{2.11}
&&D(-\nabla_x)^\top A D(\nabla_x) u=\Lambda^0 u \ \  \mbox{ in} \  \varpi (h) ,\\
\label{2.12} 
&&D(\nu(x))^\top A D(\nabla_x) u= 0  \ \  \mbox{ a.e.\,in} \ \partial \varpi (h) \,. 
\end{eqnarray}
It is a trivial consequence of the above remarks that the eigenvalues
and corresponding eigenvectors of \eqref{2.11}--\eqref{2.12} are now given by 
\bea
{\lambda_j^h} := a_h^2 \lambda_j  
= \lambda_j + O(h) \ , \ \ {u_j^h}(x)  :=  a_h^{3/2} u_j ( a_h x)  , \label{2.14}
\eea
and the relation $( {u_j^h} , u_{k}^h)_{\varpi(h)} = \delta_{j,k}$ still holds.

\BEL
\label{lem2.5} There exist  $b_1 > 0$ such that for all $k \in \bbN$, the bounds
$|u_k(x)| \leq C_k \ , \ |\nabla u_k(x)| \leq C_k 
$
hold for $x \in \varpi_0 \subset \bbR^3$ with  $|x -P^\pm | < b_1$. 
\ENL

The proof of this fact follows from standard local elliptic estimates and our smoothness
assumption of the boundary $\partial \varpi_0$ near the points $P^\pm$, see \eqref{varpi}. The lemma is proven
in detail in \cite{naruta2}, Lemma 3.1. As a corollary, we also have
\bea
|{u_k^h}(x)| \leq C_k \ , \ |\nabla {u_k^h} (x)| \leq C_k \label{2.32}
\eea
for all small enough $h> 0$, for all $x \in \varpi(h) $ with  $|x -P^\pm | < b_1$.

\subsection{Additional notation.}
\label{sec2.4}

We finish this section by fixing some  more notation. 
We write in any dimension $n$, 
 $B(x_0 ,r) = \{ x  \in \bbR^n\, : \, |x - x_0| < r \} \subset \bbR^n$
and $A(s,r) = \{ x  \in \bbR^n\, : \, s <  |x| < r \} \subset \bbR^n$ for $0 < s < r$,
and  $e^{(j)} \in \bbR^n$, $j=1,\ldots, n$, for the canonical basis vectors.

We define some smooth cut-off functions to be used at several points later. 
First, let $\chi_\theta : \bbR^3 \to [0,1]$ be a $C^\infty$-smooth function equal to 1 on the set $\theta 
\times [-1,1]$ and equal to 0 outside some neighbourhood of this set. More precisely,
we assume that $\chi_\theta$ is chosen such that 
\bea
\varpi_h \cap \{ x \, : \, \mathcal{X}^h_+(x)< 1 \} \subset \varpi(h) \ \ 
\mbox{and} \ \ 
\varpi_h \cap \{ x \, : \, \mathcal{X}^h_- (x) < 1 \} \subset \varpi(h) . \label{2.62minus}
\eea
for all  $h \in (0,\frac1{10}]$, where $\mathcal{X}^h_\pm(x)=\chi_\theta(h^{-1}(x-P^\pm))$ .
We also set for $x \in \bbR^3$
\bea
\mathcal{X}^h(x)=1-\mathcal{X}^h_-(x)-\mathcal{X}^h_+(x)  \label{4.77x},
\eea 
As a consequence,  
\bea
\varpi_h \cap \mbox{supp}\,\mathcal{X}^h  \subset \varpi(h) 
\ \ \mbox{and} \ \ 
|\nabla \cX^h | \leq {C}{h^{-1}} , \label{5.38}
\eea 
and $\cX^h$ equals 0 in both balls $B(P^\pm,ch) $ 
and 1 outside the   balls $B(P^\pm,c_1h) $ 
for some  constants $0 < c < c_1$. Also, notice that at least for small enough $h$, 
the boundary $\partial \varpi (h)$ is infinitely smooth (planar) and  \eqref{2.32} holds inside 
the domain
$\{ x \in \bbR^3 \; : \, \cX^h (x) < 1  \} .$

We also denote $\chi_\pm : \bbR^3 \to [0,1]$ some $C^\infty$-cut-off functions such that 
\bea
\chi_\pm (x) = 1 \ \mbox{for} \ x \in B( P^\pm, S ) \ \ \mbox{and} 
\ \ \chi_\pm (x) = 0 \ \mbox{for} \ x \notin B(P^\pm,  R) ; \label{4.77z}
\eea
for some  $0 < S < R < 1/2$. 
Obviously,  $\chi_{+} (y,z)= 0$ for $z \leq 0$ and $\chi_{-} (y,z)= 0$ for $z \geq 0$
and also for small enough $h$, 
\bea
 \chi_\pm \cX^h_\pm =   \cX^h_\pm  
\label{4.77zz}
\eea

\begin{figure}
\begin{center}
\includegraphics[height=8cm,width=8cm] {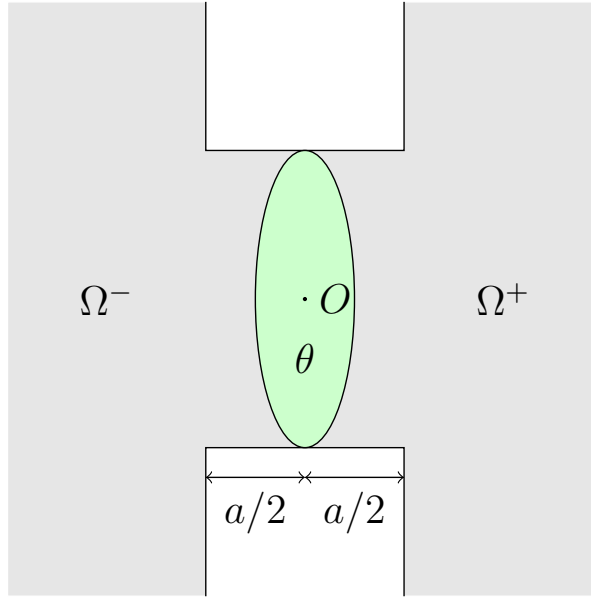} 
\end{center}
\caption{The set $\Omega$.}
\label{fig17}
\end{figure}

We define for  $a \in \{ 0,1 \}$ the domains
\bea
\Omega & := & \{ \xi \in \mathbb R^3 \, : \, |\xi_3 | > a/ 2 \} \cup
\big(  \theta \times [-a/2,a/2] \big)  \ \ 
\mbox{and } \nonumber   \\
\Omega^\pm & = & \{\xi\in \Omega:\,\, \pm \xi_3>0\} , \label{4.77w}
\eea
see Fig.\,\ref{fig17}.
The reader should keep in mind that if $x \in \varpi_h \cap {\rm supp}\, \chi_\pm$, then $\xi =h^{-1}(x-P^\pm)  \in \Omega^\mp$. 
We define the coordinate change mappings 
$\tau_\pm^h (x) = (x - P^\pm) /h$ and denote
$\tau_\pm := \tau_\pm^1$.


\section{Formal asymptotic procedure}
\label{sec5}

\subsection{The ans\"atze, the first step.} \label{sec5.1}
In this section we present 
asymptotic formulas as $h \to 0$, for the eigenvalues 
$\Lambda_k^h(\eta)$ and  for the eigenfunctions $U_k^h$, see \eqref{Lambda}, \eqref{ortho-U} and 
\eqref{5.1}, \eqref{anzatz-U}.
We start with the case that an eigenvalue  $\lambda_k$ of the limit problem is simple, see 
\eqref{lambda}. 
Let us fix such a $k$. 

As for the corresponding eigenvalue $\Lambda_k^h(\eta)$ of the model problem
\eqref{prob-2-1}--\eqref{prob-2-4}, we expect it to be a perturbation
of $\lambda_k$, and the asymptotic formula or ansatz is written as
\bea
\label{5.1}
\Lambda_k^h(\eta)&=&\lambda_k+h\Lambda_k'(\eta)+\widetilde{\Lambda}_k^h(\eta),
\eea 
where $\Lambda_k'(\eta)$ is the main correction term, to be calculated in 
Section  \ref{sec5.3}, and $\widetilde{\Lambda}_k^h(\eta)$ is a remainder, 
which will be shown to be of small order $O(h^{3/2})$ in  Section \ref{sec6}. Next, we introduce the following asymptotic formula
for the  eigenfunction $U_k^h$,
\begin{eqnarray}
\label{anzatz-U}
U_k^h(x)&=&  \mathcal{X}^h(x)\Big(u_k (a_hx)+hu_k'(a_hx ))\Big) 
\\
&+&\sum_\pm \chi_\pm(x){{V_k^\mp}}(h^{-1}(x-P^\pm) )-\sum_\pm \chi_\pm(x) \mathcal{X}^h(x)u_k(P^\pm) 
\nonumber\\
&+&\widetilde{U}_k^h (x)\ , \ \ \ x \in \varpi_h, \ \eta \in [0, 2\pi),   \nonumber
\end{eqnarray}
where 
$u_k$ is as in \eqref{ortho-u} and the notation of Section \ref{sec2.4} is used, in 
particular, 
the cut-off functions are given in \eqref{4.77x} and \eqref{4.77z},
and the function $\cX^h u_k \circ a_h$ is extended from 
$\varpi(h)$ to $\varpi_h$ as 0 by using \eqref{2.62minus}. 
The first and  second  rows of the right hand side of \eqref{anzatz-U} are called the 
outer and inner  expansions, respectively. They describe the behaviour of 
$U_k^h$ in $\varpi_h$ (at a distance of $P^\pm$)
and close to $P^\pm$, respectively; see the definitions of the cut-off functions.

To motivate the ansatz \eqref{anzatz-U} we follow the standard asymptotic
scheme described in \cite{Il} and the concordance method.
As the {\it first step} in the construction  \eqref{anzatz-U} 
we observe that according to \eqref{varpi(j,h)}, \eqref{varpi_h}, 
at some distance of the points $P^\pm$ the  domain $\varpi_h$ is just a scaling of $\varpi$, 
hence, the leading term 
$
u_k \circ a_h 
$
is expected to describe the behaviour 
of $U_k^h$ inside $\varpi_h$.

The other terms of \eqref{anzatz-U} will be considered in subsequent sections,
in particular the main correction term $u_k'= u'_k(x,\eta)$ of the outer expansion will be 
defined  in  Section \ref{sec5.3}. The behaviour of $U_k^h$ near the point $P^\pm$ is
quite  subtle, and the  main task of Section \ref{sec5.2} will be to motivate the inner
expansion, which describes  this. In particular we shall determine in Section \ref{sec5.2}
the functions ${V_k^\pm}$ and thus 
the boundary layer terms of the second row of \eqref{anzatz-U}.
The term $\widetilde U_k^h$ is again a small  remainder to be evaluated in Section 
\ref{sec6}.

\subsection{Second step: construction of  the inner expansion.} \label{sec5.2}

We next derive the expression for the inner expansion in the ansatz \eqref{anzatz-U},
\bea
\sum_\pm \chi_\pm(x){{V_k^\mp}}(h^{-1}(x-P^\pm))-\sum_\pm \chi_\pm(x) \mathcal{X}^h(x)u_k(P^\pm) .
\label{inex}
\eea 
In other words, we
aim to determine the functions ${V_k^\pm}$ and boundary layer type terms near the points 
$P^\pm$,  a task which is a priori quite unclear. The starting point is that  the above 
fixed first term $u_k \circ a_h$ does not 
satisfy the  quasiperiodicity conditions, and we thus require that 
the principal terms containing ${{V_k^\pm}}$ should compensate the discrepancy caused by that fact.
We next show how this condition will fix the functions  ${{V_k^\pm}}$, see \eqref{V} below.

By Lemma \ref{lem2.5} and the mean value theorem, 
$
u_k(a_hx) =u_k(P^\pm)+O(h) ,
$
if $x$ belongs to  some small neighbourhood $B(P^\pm ,  ch)$  of $P^\pm$.
Denoting $\xi= \tau_\pm^h(x) =h^{-1}(x- P^\pm)$, 
we have  $\xi \in \Omega^\mp$  for $x \in \varpi_h \cap {\rm supp}\, \chi_\pm$ (see \eqref{4.77w}),
and  the scaling  $L_x=h^{-2}L_\xi  = h^{-2} D(-\nabla_\xi )^\top
A D(\nabla_\xi)$ holds for the operator \eqref{4.77u}.  
So, to compensate the above discrepancy, 
we expect the functions ${{V_k^\pm}}$  to be the solutions of the problems
\begin{eqnarray}
L_\xi {{V_k^\pm}}(\xi)&=&0, \quad \xi\in \Omega^\pm ,
\label{3.21j}\\
N_\xi^{(\pm)} {{V_k^\pm}}(\xi)&=&0, \quad \xi\in \partial \Omega^\pm \setminus (\theta\times\{0\}),
\label{3.21k}
\\
{{V_k^\pm}}(\xi)&\to& u_k(P^\mp), \quad |\xi_3|\to \infty ,
\label{3.21m}
\end{eqnarray}
with additional linking boundary conditions
\begin{eqnarray}
\label{link-1}
&&V_k^{-}(\xi)=e^{i\eta}{V_k^+}(\xi) \ , \ \ 
-N_\xi^{(-)} V_k^{-}(\xi)=e^{i\eta}N_\xi^{(+)} V_k^{+}(\xi), \quad \xi\in \theta\times \{0\}.
\end{eqnarray}
Here we have denoted by $ N_\xi^{(+)} = D(\nu(\xi) )^\top
A D(\nabla_\xi)$  the  boundary operator 
associated with the domain
$\Omega^+$, and similarly for the sign "$-$". The conditions \eqref{link-1}
are related to  the quasiperiodicity conditions for $U_k^h$.

The rest of this section is devoted to looking for a solution to \eqref{3.21j}--\eqref{link-1}
and to stating some of its properties. First, note that if  a function
$V_k$ is a solution of 
\begin{eqnarray}
L_\xi V_k(\xi)&=&0, \quad \xi\in \Omega, \label{3.21}\\
N_\xi V_k(\xi)&=&0, \quad \xi\in \partial \Omega, \label{3.23} \\
V_k(\xi)&\to& u_k(P^+), \quad \xi_3<0, \quad |\xi|\to  \infty, \label{3.25} \\
V_k(\xi)&\to& e^{i\eta}u_k(P^-), \quad \xi_3>0, \quad |\xi|\to \infty,
\label{3.27}
\end{eqnarray}
where  $N_\xi= D(\nu(\xi) )^\top A D(\nabla_\xi) $ on the boundary
of $\Omega$, then the functions ${{V_k^\pm}}$ can be found as the  restrictions
\bea
{V_k^-}=V_k|_{\Omega^-}, \quad {V_k^+}=e^{-i\eta}V_k|_{\Omega^+}.
\label{5.37c}
\eea

\BED
\label{def5.1} 
For $j=1,2,3$,
we denote by  $X_j : \Omega \to \bbR^3 $ the three functions which are
the solutions of the problems
\begin{eqnarray}
L_\xi X_j(\xi)&=&0, \quad \xi\in \Omega, \label{5.22} \\
N_\xi X_j(\xi)&=&0, \quad \xi\in \partial \Omega, \label{5.22y}
\end{eqnarray} 
with asymptotics $X_j(\xi)\to \pm e^{(j)}$ when $\pm\xi_3>0$ and $|\xi|\to \infty$.
\END

We shall soon prove the existence of  such functions 
$X_j$. Taking that for granted we define the function $V_k$ by 
\bea
\label{5.37}
V_k(\xi)=\sum\limits_{j=1}^3 \big(a^{(k)}_j(\eta)e^{(j)}+b_j^{(k)}(\eta)X_j(\xi)\big)
\eea
where the coefficients are to be chosen such that \eqref{3.21}--\eqref{3.27} are
satisfied. Since the functions $X_j$ do not depend on $\eta$ or $k$, 
it is plain that $V_k$ depends on $\eta$ and $k$ only via the coefficients, which are determined
by the  six equations
\begin{equation}
\label{sys-a-b}
a^{(k)}(\eta)-b^{(k)}(\eta)=u_k(P^+), \quad a^{(k)}(\eta)+b^{(k)}(\eta)=e^{i\eta}u_k(P^-)
\end{equation} 
with  $a^{(k)}=(a^{(k)}_1,a^{(k)}_2,a^{(k)}_3)^\top$ and 
$b^{(k)}=(b_1^{(k)},b_2^{(k)},b_3^{(k)})^\top$. Solving the system \eqref{sys-a-b}
yields
\bea
\label{a-eta} \ \ \ \ 
a^{(k)}(\eta)&=&\frac{1}{2}(u_k(P^+)+e^{i\eta}u_k(P^-))\ , \  
b^{(k)}(\eta) = \frac{1}{2}(-u_k(P^+)+e^{i\eta}u_k(P^-)).
\eea

We now turn  to the existence and some properties of the functions $X_j$.

\BEL
\label{lemma5.2} For every $j=1,2,3$, there exists a unique solution to the problem \eqref{5.22}--\eqref{5.22y}, and it can be written in the form
\bea
\label{Xj}
X_j & = & \pm \big(1 - \chi_\theta \big)   e^{(j)} 
+ \big(1 - \chi_\theta \big)  \sum_{l=1}^3M^\pm_{lj}T_l^\pm \circ \tau_\pm 
 + \widetilde{T}_j 
\ \ \mbox{in} \ \Omega^\pm ,
\eea 
where 
 $ \tau_\pm(\xi)  = \xi - P^\pm$ and $M^\pm_{lj} \in \bbR$ are unique coefficients 
and   $T_l^\pm$ denote
the  Poisson kernels  for the operator $L_\xi$ in the domains 
$\{ \xi \in \bbR^3 \, : \, \pm \xi_3 > 0 \}$, respectively (see \eqref{5.24}--\eqref{5.24xy}). 
The remainders $\widetilde{T}_j$ satisfy for $\xi \in \bbR^3$ 
the estimates 
\begin{equation}
\label{remainder-X}
|\widetilde{T}_j(\xi)|\leq \frac{C}{1+ |\xi|^2},
\quad |\nabla_\xi \widetilde{T}_j(\xi)|\leq \frac{C}{ 1+ |\xi|^3} . 
\end{equation}

\ENL


Proof.  We define the  Poisson kernels $T_l^\pm$ as 
solutions of the half-space  problems
\begin{eqnarray}
L_\xi T_l^\pm (\xi) &=&0, \quad \xi\in \bbR^3, \  \pm\xi_3> 0 , \label{5.24} 
\\
N_\xi T_l^\pm (\xi) &=&\delta_{(\xi_1,\xi_2)} e^{(l)}, \quad \pm\xi_3= 0.
\label{5.24xy}
\end{eqnarray}

Here, the differential operators $L_\xi$ and $N_\xi$ are understood as distributional 
derivatives (notice the direction of the normal vector in \eqref{5.24xy}, e.g. it is $-e^{(3)}$ in the case of the sign $"+"$)
and $\delta_{(\xi_1,\xi_2)}$ stands for the 2-dimensional 
Dirac measure of the point $(0,0)$ in the plane $\{ \xi \in \bbR^3 \, : \, \xi_3 = 0\}$.
The solutions of the problem \eqref{5.24} are known to be homogeneous functions of order 
$-1$, in other words, the function $ T_l^\pm$ has a
singularity at   $0$, and  for 
all $\xi \in \bbR^3$, $ \pm\xi_3>0$, $0 < r \in \bbR$,  we have
\bea
& & T_l^\pm ( r \xi) = \frac1r T_l^\pm(\xi) \ \ , \ \  | T_l^\pm( \xi) |  
\leq  \frac{C}{|\xi|}
 \ \ , \ \  | \nabla_\xi T_l^\pm( \xi) |  
\leq  \frac{C}{|\xi|^2} . 
   \label{5.24a}
\eea

Let us define for a moment a smooth 
function in   $\Omega$
by  $\Theta^{(j)} (\xi) =  \pm \big(1 - \chi_\theta(\xi)\big)e^{(j)}$,
where $ \xi \in \Omega^\pm$ and $\chi_\theta$ is as in \eqref{2.62minus}.
We now look for  a function $X_j : \Omega \to \bbR$ in the form 
$
X_j =  \Theta^{(j)} 
+\widetilde{X}_j 
$,
where  the function $\widetilde{X}_j$ satisfies the problem
\begin{eqnarray*}
L_\xi \widetilde{X}_j(\xi)&=& - L_\xi \Theta^{(j)}(\xi), \quad \xi\in \Omega, \\
N_\xi \widetilde{X}_j(\xi)&=& -  N_\xi \Theta^{(j)}(\xi), \quad \xi\in \partial \Omega ,
\end{eqnarray*} 
and the functions on the right are smooth and have compact supports. According
to the general results from \cite{KoOl}, such a function $\widetilde{X}_j$ exists and
is unique in the function space  
\begin{equation}
\label{energy-X}
\big\{ v \in H_{\rm loc}^1 ( \Omega) \, : \, 
\|\nabla_\xi v ; L^2(\Omega)\|+\|(1+|\xi| )^{-1}v; 
L^2(\Omega)\|< \infty \big\} . 
\end{equation}
Using the results of \cite{na262} and \cite{MaPl2}, the solution 
$\widetilde{X}_j$ is, say in a bounded ball $B(0,2)$, at least 
$C^2$-smooth bounded
function, and outside $B(0,2)$, a linear  combination of the 
Poisson kernels plus a perturbation which is small at the infinity.
More precisely, we can write for  some coefficients $M_{lj}^\pm$,
\begin{eqnarray}
& & \widetilde{X}_j (\xi)= (1 - \chi_\theta(\xi)) \sum_{l=1}^3 M_{lj}^\pm
T_l^\pm (\xi - P^\pm ) + 
\widetilde{T}_j(\xi) \ , \ \xi \in \Omega^\pm , 
\label{5.37u} 
\end{eqnarray}
where the functions $ (1 - \chi_\theta) T_l^\pm (\xi - P^\pm)$
are well defined in $\Omega^\pm$ and 
the perturbation $\widetilde{T}_j$ satisfies \eqref{remainder-X};
moreover, $\widetilde T_j |_{B(0,r)} \in L^2(\Omega^\pm)$
for any fixed radius $r >0$.
More information on the 
coefficients $ M_{lj}^\pm$ will be obtained in Lemma \ref{lemma5.3} below.
 \ \ $\Box$

\bigskip

We now define $V_k$ by putting  \eqref{Xj}  into \eqref{5.37} and using
\eqref{a-eta}. We also denote, keeping in mind that
$b^{(k)}$ and thus the consequent functions depend also on $\eta$, 
\begin{eqnarray}
\label{W-}
{W_k^\pm}(\xi)&=&  
\big(1 - \chi_\theta \circ \tau_\mp(\xi) \big) \sum_{j,l=1}^3b^{(k)}_jM_{lj}^\pm T_l^\pm(\xi)\ 
, \ \ \tau_\mp(\xi) \in \Omega^\pm  ,\\
\widetilde{W}_k^\pm  (\xi)  &=&   \sum_{j=1}^3b^{(k)}_j
\big( \mp \chi_\theta (\xi) e^{(j)} + \widetilde{T}_j(\xi)  \big) 
, \ \ \xi \in \Omega^\pm , \label{wideW} 
\end{eqnarray}
where $ \tau_\pm(\xi)  = \xi - P^\pm$.
We remark that due to \eqref{5.24a}, \eqref{W-}, 
we have $ {W_k^\pm} \in L^2(\Omega^\pm) $. By the remarks on
$\widetilde  T_j $ 
in the proof of Lemma \ref{lem5.2}, 
$\widetilde{W}_k^\pm \big\vert_{B(0,r)}  \in L^2(\Omega^\pm) $
for any  constant radius $r >0$,  and moreover 
\bea
\Vert {W_k^\pm} ;L^2(\Omega^\pm) \Vert \leq C_k \ \ , \ \ 
\Vert\widetilde{W}_k^\pm \big\vert_{B(0,r)} ;L^2(\Omega^\pm) \Vert \leq C_k .
 \label{6.3dx}
\eea
We finally find  the inner expansion by 
taking the restrictions \eqref{5.37c}, which leads to
\begin{equation}
\label{V}
{V_k^\pm}(\xi)=   u_k(P^\mp)+ 
{W_k^\pm} \circ \tau_\pm(\xi)  
+\widetilde{W}_k^\pm  (\xi) 
, \quad \xi \in \Omega^\pm .  
\end{equation}
From  \eqref{W-}, \eqref{5.24a} we see that there exists 
 some  constant $M >2$, which we fix now, such that 
for all $\xi \in \bbR^3$ with $|\xi| \geq M/2$ 
\bea
& & 
{W_k^\pm} ( r \xi) = \frac1r {W_k^\pm}  (\xi) 
 \  ,  \label{5.24af} 
\eea
Also, \eqref{5.24a}, \eqref{remainder-X} imply for all $\xi \in \bbR^3$ 
\bea
& & 
|  {W_k^\pm}  (\xi)   | \leq \frac{C_k}{1 +  |\xi |} 
  \ \ , \ \  
 | \nabla_\xi {W_k^\pm}  (\xi) 
 | \leq \frac{C_k}{ 1+ |\xi |^2} 
\ ,  \label{5.24afNEW} \\ 
& & 
|\widetilde{W}_k^\pm  (\xi)
|\leq \frac{C_k}{1 + |\xi |^2} \ , \ \ 
|\nabla_\xi \widetilde{W}_k^\pm  (\xi )
|\leq \frac{C_k}{1 + |\xi|^3}.
\label{5.82}
\eea

We complete the present study of the inner expansion by the following remarks. 
\BEL
\label{lemma5.3}
The matrices ${\bf M}^\pm=(M_{jl}^\pm)_{j,l=1}^3$ defined by the coefficients in \eqref{Xj},  \eqref{5.37u} (which do not depend on $\eta$ or
$k$)  have the properties 
\bea 
\label{5.20}
{\bf M}^+=-{\bf M}^{-} 
\ \ , \ \ 
{\bf M}^\pm = \big( {\bf M}^\pm\big)^\top , 
\eea
and, moreover,  $\bfM^+$ is positive definite.
\ENL

Proof. Let us use the (generalized) Green formula (see \cite{Lad})
for $X_l$ and $e^{(j)}$ in the domain $B (0,\varrho)\cap \Omega$, where
$\varrho > 0$ large enough, in particular such that $\chi_\theta (\xi) = 0$ 
for $|\xi| \geq \varrho - 1/2$. So, by \eqref{5.22}, \eqref{5.22y} and Fig.\,\ref{fig17}.,
\begin{equation*}
0=(L_\xi X_l, e^{(j)})_{B (0,\varrho)\cap \Omega}= (N_\xi X_l, e^{(j)})_{S_\varrho^+}+(N_\xi X_l, e^{(j)})_{S_\varrho^-} ,
\end{equation*}
where $S_\varrho^\pm=\{\xi: |\xi|=\varrho, \frac{a}{2}<\pm\xi_3\}$. Taking into account the representation \eqref{Xj}, we get
\bea
0 & = & (N_\xi X_l, e^{(j)})_{S_\varrho^+}+(N_\xi X_l, e^{(j)})_{S_\varrho^-}=
\sum_{\pm}\sum_{m=1}^3 M_{ml}^\pm (N_\xi T_m^\pm \circ \tau_\pm, e^{(j)})_{S_\varrho^\pm}
\roweq
- \sum_{\pm}\sum_{m=1}^3 M_{ml}^\pm(N_\xi T_m^\pm\circ \tau_\pm, e^{(j)})_{\{\pm\xi_3=a/2\}}
\roweq
- \sum_{\pm}\sum_{m=1}^3 M_{ml}^\pm(N_\xi T_m^\pm, e^{(j)})_{\{\pm\xi_3=0\}}
= 
-M_{jl}^+-M_{jl}^-\, , \nonumber
\eea
where  we also used  \eqref{5.24xy} and it was necessary to change the sign due to 
the direction of the normal vector. 
Notice that the contribution of the term $\widetilde T_m$ must vanish
as a consequence of the second inequality \eqref{remainder-X}, since 
$(N_\xi \widetilde T_m, e^{(j)})_{S_\varrho^\pm} = O(\varrho^{-1})$ for large $\varrho$,
whereas the terms $(N_\xi  T_m^\pm, e^{(j)})_{S_\varrho^\pm}$  are constant
with respect to $\varrho$. We have proven the first identity in \eqref{5.20}.

To prove the second identity of \eqref{5.20} we use the Green formula for the functions $X_j$
and $X_l$:
\bea
0 & = & (L_\xi X_j, X_l)_{B (0,\varrho)\cap \Omega}= (N_\xi X_j, X_l)_{S_\varrho^+}
\rowpl
(N_\xi X_j, X_l)_{S_\varrho^-} -(X_j, N_\xi X_l)_{S_\varrho^+}-(X_j, N_\xi X_l)_{S_\varrho^-}. \label{3.35z}
\eea
Since $N_\xi T_m^\pm (\xi)  = O(\varrho^{-2})$ and $T_m^\pm (\xi) = O(\varrho^{-1})$ 
for $\xi \in S_\varrho$,  the scalar products 
$(N_\xi T_m^\pm, T_l^\pm)_{S_\varrho^\pm}$ tend to zero as $\varrho\to +\infty$. 
Thus, by \eqref{3.35z} and the same observations as in the first case, 
\bea
0 & =& \lim_{\varrho\to+\infty}\sum_{m=1}^3\Big(M_{mj}^+(N_\xi T_m^+, e^{(l)})_{S_\varrho^+}+M_{mj}^-(N_\xi T_m^-, -e^{(l)})_{S_\varrho^-} 
\rowmi
M_{ml}^+(e^{(j)}, N_\xi T_m^+)_{S_\varrho^+}-M_{ml}^-(-e^{(j)}, N_\xi T_m^-)_{S_\varrho^-}\Big)=2M_{lj}^--2M_{jl}^+.
\eea

Finally, for the positive definiteness of $\bfM^+$ we 
consider a vector $v=(v_1,v_2,v_3)^\top \in \bbR^3$ and the function $w=\sum_{j=1}^3v_jX_j$.
Then, by the positivity of the operator $L_\xi$, 
\bea
0  & \leq & (L_\xi w,w)_ {\Omega^+} =
(N_\xi w, w)_{\partial \Omega^+}= -(N_\xi w, w)_{S_\varrho^+}
\roweq
-\sum_{j,l=1}^3v_j \overline{v_l}\int\limits_{S_\varrho^+} X_l(\xi) N_\xi X_j(\xi)\, ds_\xi\,. \label{3.35zd}
\eea
It remains to note as above that 
$$
\lim_{\varrho\to+\infty}\int\limits_{S_\varrho^+} X_l(\xi) N_\xi X_j(\xi)=-M^+_{jl}\,,
$$
hence,  $v^\top {\bf M}^+ \overline{v} = 
(L_\xi w,w)_ {\Omega^+}$. The claim follows from \eqref{3.35zd}. \ \ $\Box$

\subsection{Third step: construction of the outer expansion and the 
leading correction term for the eigenvalue.}
\label{sec5.3}
We next look for the correction term $u_k'$ of the outer expansion.
Let us write for a moment $E(x,h) :=  u_k(a_hx)+ h 
\tilde u_k (a_h x,\eta)$ for some $\tilde u_k \in H^1(\varpi_0)$. 
Since $L_x u_k (x) = \lambda_k u_k(x)$,  we can calculate
\bea
& &  L_x E(x,h) 
=  
\lambda_k u_k(a_hx) + 2 a h \lambda_k u_k(a_hx)+ 
h L_x \tilde u_k(a_h x,\eta)  + O(h^2)
\nonumber 
\eea
and on the other hand, assuming \eqref{5.1}, 
\bea
& & \Lambda_k^h(\eta)   E(x,h)
= 
\lambda_k u_k(a_hx)  +  h \Lambda_k'(\eta) u_k(a_hx) +  
\lambda_k h \tilde u_k(a_h x) + O(h^{3/2}) 
\nonumber 
\eea
Hence, the assumption that $L_x E(x,h) = \Lambda_k^h E(x,h) $ is  correct  up to terms of order $h^{3/2}$,
leads to the following problem   to find $\tilde u_k$ in $\varpi_0$:
\begin{eqnarray}
&&(L_x-\lambda_k) \tilde u_k(x,\eta)=(-2a\lambda_k+\Lambda_k'(\eta))u_k(x),\quad x\in\varpi_0,
\label{5.92a}  \\
&&{N_x} \tilde u_k(x,\eta)=0, \quad x\in \partial \varpi_0,   \label{5.92b}   \\
&&\tilde u_k(x,\eta)\sim  {W_k^\pm}(x-P^\mp), \quad x\to P^\mp.  \label{5.92c}
\end{eqnarray}
The asymptotic condition \eqref{5.92c} comes by looking at  the inner expansion \eqref{inex}
near  the points $P^\pm$.  We write   $\tilde u_k$ as 
\beq
\tilde u_k(x,\eta)=\sum_\pm\chi_\pm(x) 
{W_k^\mp}(x-P^\pm)
 +  u_k'(x,\eta), \label{5.83a}
\eeq
where  $\chi_\pm$ are  as in \eqref{4.77z},
and state the following problem for $ u_k'$:
\bea
(L_x-\lambda_k) u_k'(x,\eta) 
&=& F(x)
, \quad x\in\varpi_0,
\label{prob-U'-1}
\\
{N_x}  u_k'(x,\eta) 
& = & G(x) 
, \quad x\in \partial\varpi_0.  \label{prob-U'-2}
\eea
Here we denote (cf. \eqref{4.77z})
\bea
F & = &  \big(-2a\lambda_k+\Lambda_k'(\eta)\big) u_k
- (L_x-\lambda_k)
\sum_{\pm}\chi_{\pm}
{\cW_k^\mp} \circ \tau_\pm  \  \ \mbox{in} \ \varpi_0, 
\label{Fdef}
\eea
\begin{equation}
G(x) = \left\{
\begin{array}{cl}
-\sum_\pm {N_x}(\chi_\pm(x)  
{\cW_k^\mp}(x-P^\pm) \ , & x \notin  \partial \varpi_0 \cap B(P^\pm , S/2) 
\\
0  \ , & x \in  \partial \varpi_0 \cap B(P^\pm , S/2) ,
\end{array}
\right. \label{Gdef}
\end{equation}
where $\cW_k^\mp$ are defined as $W_k^\mp$ in \eqref{W-} but
without the cut-off function $\chi_\theta$:
\bea
\label{cW}
{\cW_k^-}(\xi)&=&    \sum_{j,l=1}^3b^{(k)}_jM_{lj}^- T_l^-(\xi)\ 
, \ \ \xi \in \Omega^-  ,
\eea
and similarly for $W_k^+$.
Note that  $\chi_\pm (x) =1$ for $x$ close to  $P^\pm$.
Hence,  \eqref{5.24},  \eqref{5.24a}, \eqref{W-}
imply  that  the functions 
$
L_x(\chi_\pm  W_k^\mp \circ \tau_\pm ) 
$ 
are equal to 0 near the points $P^\pm$, and thus $F\in L^q(\varpi_0)$ for every $q<3$.
The function  $G$ is smooth and thus bounded, by similar arguments and \eqref{5.24xy};
in fact, $G$ equals the expression on the first row of \eqref{Gdef} everywhere
except at $P^\pm$, where the latter has the Dirac measure singularity.

\BEL
\label{lem3.5}
If 
\begin{equation}
\label{L'}
\Lambda_k'(\eta)=2a\lambda_k+2b^{(k)}(\eta)^\top {\bf M}^+ 
\overline{b^{(k)}(\eta)}
\end{equation}
then the problem \eqref{prob-U'-1}-\eqref{prob-U'-2} has a unique solution in $H^1(\varpi_0)^3$.
\ENL

{\it Proof.} 
The problem \eqref{prob-U'-1}-\eqref{prob-U'-2} can be rewritten in the weak formulation as 
\begin{equation}
\label{weak-U'}
(AD(\nabla_x)u_k',D(\nabla_x)v)_{\varpi_0}-(G,v)_{\partial \varpi_0}-\lambda_k(u_k',v)_{\varpi_0}=
(F,v)_{\varpi_0},\quad v\in H^1(\varpi_0)^3.
\end{equation}
Equation \eqref{weak-U'} is equivalent to 
\begin{equation}
\langle u_k',v \rangle_{\cH^0}
-(1+\lambda_k)\langle \mathcal{B}^0u_k',
v \rangle_{\cH^0}=
(G,v)_{\partial \varpi_0}+(F,v)_{\varpi_0} \ , \quad v\in\mathcal{H}^0,
\nonumber
\end{equation}
where $\cB^0$ and $\cH^0$ were defined  below \eqref{prob-3-2}. 
We use the fact that functional defined by the formula
$
\cF (v)=(G,v)_{\partial \varpi_0}+(F,v)_{\varpi_0}
$ 
is linear and continuous on $H^1(\varpi)^3$. 
The problem \eqref{prob-U'-1}-\eqref{prob-U'-2} can thus be rewritten as the equation
$
\mathcal{B}^0 u_k' - \mu_k u_k'=-\mu_k \mathcal{F},
$
where $\mu_k=(1+ \lambda_k)^{-1} $.
According to the  Fredholm alternative,  this equation has a solution if and only if
the right-hand side is orthogonal to $u_k$ (solution of the homogeneous
problem, Section \ref{sec2.3}). So we get the  solvability condition
\bea
0&=&(G,u_k)_{\partial \varpi_0}+(F,u_k)_{\varpi_0} . \label{solva}
\eea
Here we use the Green formula and take into account \eqref{prob-3-1} and \eqref{prob-3-2}:
\begin{eqnarray*}
& & (F,u_k)_{\varpi_0}=
-2a\lambda_k+\Lambda_k'(\eta)
-\sum_\pm \int\limits_{\varpi_0}(L_x-\lambda_k)(\chi_\pm{\cW_k^\mp} \circ \tau_\pm )
 \overline{u}_kdx
\roweq
-2a\lambda_k+\Lambda_k'(\eta)
+ \sum_\pm \int\limits_{\varpi_0}\lambda_k \chi_\pm{\cW_k^\mp} \circ \tau_\pm   \overline{u}_k dx
- 
\sum_\pm \int\limits_{\varpi_0}\chi_\pm {\cW_k^\mp} \circ \tau_\pm L_x\overline{u}_kdx
\rowmi
\sum_\pm \int\limits_{\partial \varpi_0}N_x^0 ( \chi_\pm{\cW_k^\mp} \circ \tau_\pm ) 
\overline{u}_kds_x
+ 
\sum_\pm \int\limits_{\partial \varpi_0}\chi_\pm{\cW_k^\mp} \circ \tau_\pm N_x^0
\overline{u}_kds_x
\roweq
-2a\lambda_k+\Lambda_k'(\eta)
- \sum_\pm \int\limits_{\partial \varpi_0}N_x^0 ( \chi_\pm{\cW_k^\mp} \circ \tau_\pm ) 
\overline{u}_kds_x .
\end{eqnarray*} 
By the remarks just before Lemma \ref{lem3.5},  the last integral equals 
$(G,u_k)_{\partial \varpi_0}$ except that $N_x^0 ( \chi_\pm{\cW_k^\mp} \circ \tau_\pm ) $
contains Dirac measures at $P^\pm$ which $G$ does not contain. This remark, 
\eqref{5.24xy}, \eqref{W-},   and \eqref{solva} yield
\bea
& & -2a\lambda_k+\Lambda_k'(\eta)
= 
b(\eta)^\top {\bf M}^- u_k(P^+)+e^{-i\eta} b(\eta)^\top {\bf M}^+ u_k(P^-)=
2b(\eta)^\top {\bf M}^+ \overline{b(\eta)} .  \ \  \Box \nonumber 
\eea

\BER
\label{rem3.6}
The facts that  $F \in L^q(\varpi_0)$ for  every $q<3$, the function  $G$ is smooth,
and the boundary of $\varpi_0$ is smooth near points $P^\pm$, imply that $u_k' \in W^2_q(\varpi_0)$, by standard elliptic estimates.
Due to the  Sobolev embedding theorem and \eqref{5.83a} we get  for all $p<\infty$
that $\nabla_x u_k'\in L^p(\varpi_0) $ and $ u_k'\in L^\infty(\varpi_0) $
with the corresponding norm bounds independent of $\eta$.
\ENR

\subsection{Comments on multiple eigenvalues.}
\medskip
Let us consider the case of a multiple  eigenvalue $\lambda_k$, see \eqref{lambda}.
Assume that its multiplicity is $m \geq 2$ and let
$u_{k+j}$, where $j =0,\ldots, m-1$, be the  corresponding eigenfunctions satisfying the 
orthogonality and normalization condition \eqref{ortho-u}. Now the 
principal term of the asymptotic of $U_{k+l}^h$, where $l =0,\ldots, m-1$, 
is a linear combination
$
\sum_{j=0}^{m-1} \alpha_{l,j}  u_{k+j} \circ a_h =: u_{k,l}\circ a_h .
$
Repeating the procedure of the concordance method of asymptotic expansions 
in this case we again get the  problem  \eqref{prob-U'-1}-\eqref{prob-U'-2} for the main
 correction terms $u_{k+l}'$ and $\Lambda_{k+l}'(\eta)$.
Now it turns out that we have $m$ solvability conditions which are equivalent to 
a system of $m$ linear equations 
$$
(-2a\lambda_k+\Lambda_{k+l}'(\eta))\alpha_{l,j}=
2\sum_{q=0}^{m-1} \alpha_{l,q} b^{(k+q)}(\eta)^\top 
{\bf M}^+ \overline{b^{(k+j)}(\eta)},
$$
where $j =0,\ldots, m-1$, while ${\bf M}^+$ is as in Section \ref{sec5.2} and 
$b^{(k+q)}(\eta)=\frac{1}{2}(u_{k+q}(P^-)-e^{-i\eta}u_{k+q}(P^+))$ (see \eqref{a-eta}). 
As a consequence, the expression $-2a\lambda_k+\Lambda_{k+l}'(\eta)$ is one 
of the $m$ eigenvalues (multiplicities counted) of the $m\times m$ matrix
$$
{\bf B}(\eta)= \Big( b^{(k+q)}(\eta)^\top {\bf M}^+ \overline{b^{(k+j)}(\eta)} \Big)_{q,j=1}^m\ \ ,
$$
and the coefficient sequence $(\alpha_{l,j})_{j=1}^m $ is found as the 
eigenvector of $\bfB(\eta)$ corresponding to the eigenvalue
$-2a\lambda_k+\Lambda_{k+l}'(\eta)$. 
It is clear that the rank of the matrix ${\bf B}$ does not exceed 3.

\medskip

\subsection{The case $\lambda=0$}
The sequence \eqref{lambda} begins with six eigenvalues equal to 0 corresponding to the rigid
motions
$$
u_j(x)=\beta_j e^{(j)} \ \mbox{for} \  j=1,2,3, \ \ \ \ 
u_j(x)=\beta_j x\times e^{(j-3)} \ \mbox{for} \ j=4,5,6,
$$ 
where $\beta_j$ are the normalization multipliers, i.e. $\beta_1=\beta_2=\beta_3=|\varpi_0|^{-1/2}$,
$\beta_4=(J_2+J_3)^{-1/2}$, $\beta_5=(J_1+J_3)^{-1/2}$, $\beta_6=(J_2+J_1)^{-1/2}$. Here $J_k$
is the moment of inertia of the body $\varpi_0$ around the axis $e^{(k)}$.  
Calculating the vectors $b^{(k)}(\eta)$ and suppressing the inessential
index $l$ of $\alpha_{l,j}$, we get
\begin{eqnarray*}
b^{(k)}(\eta)&=&2^{-1}|\varpi_0|^{-1/2}(1-e^{-i\eta})e^{(k)}=:\alpha_1 e^{(k)},\quad k=1,2,3;\\
b^{(4)}(\eta)&=&0;\\
b^{(5)}(\eta)&=&-4^{-1}|J_1+J_3|^{-1/2}(1+e^{-i\eta})e^{(3)}=:\alpha_5 e^{(3)};\\
b^{(6)}(\eta)&=&4^{-1}|J_1+J_2|^{-1/2}(1+e^{-i\eta})e^{(2)}=:\alpha_6 e^{(2)}.
\end{eqnarray*}
Thus, using the symmetry of matrix ${\bfM}^+$,
\begin{equation*}
\bfB(\eta)=
\begin{pmatrix}
|\alpha_1|^2M^+_{11}& |\alpha_1|^2M^+_{12} &|\alpha_1|^2M^+_{13} & 0 
& \alpha_1\overline{\alpha_5}M^+_{13} & \alpha_1\overline{\alpha_6}M^+_{12}\\
|\alpha_1|^2M^+_{12}& |\alpha_1|^2M^+_{22} &|\alpha_1|^2M^+_{23} & 0 
& \alpha_1\overline{\alpha_5}M^+_{23} & \alpha_1\overline{\alpha_6}M^+_{22}\\
|\alpha_1|^2M^+_{13}& |\alpha_1|^2M^+_{23} &|\alpha_1|^2M^+_{33} & 0 
& \alpha_1\overline{\alpha_5}M^+_{33} & \alpha_1\overline{\alpha_6}M^+_{23}\\
0&0&0&0&0&0 \\
\overline{\alpha_1}\alpha_5 M^+_{13}& \overline{\alpha_1}\alpha_5 M^+_{23} &\overline{\alpha_1}\alpha_5 M^+_{33} & 0 
& |\alpha_5|^2M^+_{33} & \alpha_5\overline{\alpha_6}M^+_{23}\\
\overline{\alpha_1}\alpha_6 M^+_{12}& \overline{\alpha_1}\alpha_6 M^+_{22} &\overline{\alpha_1}\alpha_6 M^+_{23} & 0 
& \alpha_6\overline{\alpha_5}M^+_{23} & |\alpha_6|^2M^+_{22}
\end{pmatrix}.
\end{equation*}
This can be rewritten in shorter form
\begin{equation*}
\bfB(\eta)=
\begin{pmatrix}
 & &  & 0  &  & \\
& |\alpha_1|^2{\bfM}^+ & & 0 & \alpha_1\overline{\alpha_5}{\bfM}^+_{3} & 
\alpha_1\overline{\alpha_6}{\bfM}^+_{2}\\
&  &  & 0  &  &\\
0&0&0&0&0&0 \\
& \overline{\alpha_1}\alpha_5 ({\bfM}^+_{3})^\top & & 0 
& |\alpha_5|^2M^+_{33} & \alpha_5\overline{\alpha_6}M^+_{23}\\
& \overline{\alpha_1}\alpha_6 ({\bfM}^+_{2})^\top & & 0 
& \alpha_6\overline{\alpha_5}M^+_{32} & |\alpha_6|^2M^+_{22}
\end{pmatrix} ,
\end{equation*}
where ${\bfM}^+_{k}$ is the $k$th column of the matrix ${\bfM}^+$.
If $e^{(j)}$, $1\leq j\leq 6$, are 
the standard basis vectors in $\mathbb{R}^6$ it is easy to see that the vectors
$$
e^{(4)}, \overline{\alpha_5} e^{(3)}-\overline{\alpha_1}e^{(5)},  \overline{\alpha_6}
e^{(2)}-\overline{\alpha_1} e^{(6)}
$$
are in the kernel of the matrix $\bfB(\eta)$. So, the corresponding linear combinations of the functions
$u_j$ form  the ansatz of the eigenfunctions of the perturbed problem, but the asymptotic corrections 
of order $h$ for the eigenvalues are 0. To find the  other three corrections of 
order $h$ we must find the nonzero eigenvalues of the matrix $\bfB(\eta)$. To do this we have to 
calculate the characteristic function $\phi(t)=\det (\bfB(\eta)-t {\bfE}_6)$, where $\bfE_n$ denotes the unit matrix of dimension $n \times n$. One can write

\bea
\phi(t) &=&
\det\begin{pmatrix}
 & &  & 0  & 0 & 0 \\
& |\alpha_1|^2{\bf M}^+-t{\bfE}_3 & & 0 & 0 & 
\overline{\alpha_1^{-1}\alpha_6}t\\
&  &  & 0  & \overline{\alpha_1^{-1}\alpha_5}t & 0\\
0&0&0&-t&0&0 \\
& \overline{\alpha_1}\alpha_5 ({\bfM}^+_{3})^\top & & 0 
& -t & 0\\
& \overline{\alpha_1}\alpha_6 ({\bfM}^+_{2})^\top & & 0 
& 0 & -t
\end{pmatrix} 
\roweq
\det\begin{pmatrix}
 & &  & 0  & 0 & 0 \\
& \widetilde {\bf M} -t{\bfE}_3 & & 0 & 0 & 0\\
&  &  & 0  & 0 & 0\\
0&0&0&-t&0&0 \\
& \overline{\alpha_1}\alpha_5 ({\bfM}^+_{3})^\top & & 0 
& -t & 0\\
& \overline{\alpha_1}\alpha_6 ({\bfM}^+_{2})^\top & & 0 
& 0 & -t
\end{pmatrix} = -t^3\det(\widetilde {\bf M}-t{\bfE}_3), \nonumber
\eea 
where
\bea 
\widetilde {\bf M} & =& 
\begin{pmatrix}
|\alpha_1|^2M^+_{11}& |\alpha_1|^2M^+_{12} &|\alpha_1|^2M^+_{13} \\
(|\alpha_1|^2+|\alpha_6|^2)M^+_{12}& (|\alpha_1|^2+|\alpha_6|^2)M^+_{22} &(|\alpha_1|^2+|\alpha_6|^2)M^+_{23} \\
(|\alpha_1|^2+|\alpha_5|^2)M^+_{13}& (|\alpha_1|^2+|\alpha_5|^2)M^+_{23} &(|\alpha_1|^2+|\alpha_5|^2)M^+_{33} 
\end{pmatrix} 
\roweq
\diag(|\alpha_1|^2, |\alpha_1|^2+|\alpha_6|^2, |\alpha_1|^2+|\alpha_5|^2) {\bf M}^+ \nonumber
\eea 
and $\diag(a,b,c)$ is a diagonal matrix with elements $a$, $b$ and $c$. 

So in the case $\lambda=0$ there are  six main asymptotic corrections  $\Lambda_j'(\eta)$: for 
$j=1,2,3$ they are  the eigenvalues 
of the matrix $\widetilde {\bf M}$, and for $j=4,5,6$ they are just $0$.

\section{Main result: position of spectral gaps.}
\label{sec6}

\subsection{Main theorem}
\label{sec6.1}
We state our main result on the position of the spectral bands 
for the linear elasticity problem in the domain $\Pi_h$. Recall from \eqref{eq1.1} that the bands consist of  the eigenvalues $\Lambda_k^h (\eta)$.

\BET
\label{th6.2}
For every $k$ there exists a constant $C_k >0$ such that for all  $h>0$, $\eta \in 
[0, 2 \pi)$,
\bea
& & |\Lambda_k^h(\eta) - ( \lambda_k + h \Lambda_k'(\eta) )  | \leq C_k h^{3/2}  
\ \ , \ \ \ \Vert U_k^h - \cU_{h,k} ; \cH^{h,\eta} \Vert  \leq C_k h^{3/2} .
\label{6.12}
\eea
Here $\lambda_k + h \Lambda_k'(\eta)$ is the approximate eigenvalue from
\eqref{5.1}, $\lambda_k$ is the $k$th eigenvalue of the limit problem, see \eqref{lambda},
and $\Lambda_k'(\eta)$ is determined in Lemma \ref{lem3.5} of Section \ref{sec5.3},
\bea
\Lambda_k'(\eta) = 2 a \lambda_k  + \big(A_k + e^{i\eta} B_k\big)^\top
\bfM^+ \big(A_k + e^{-i\eta} B_k\big) , \label{1.98}
\eea
where the column vectors $A_k$, $B_k \in \bbR^3$ and the positive definite matrix 
$\bfM^+ \in \bbR^{3 \times3}$ do not depend on $h$ or $\eta$. 
Moreover, 
\begin{eqnarray}
\cU_{h,k}&=&  \mathcal{X}^h\big(u_k \circ a_h +hu_k'\circ a_h \big) 
+ 
\sum_\pm \chi_\pm {V_k^\mp} \circ \tau_\pm^h - 
\sum_\pm \chi_\pm  \mathcal{X}^h u_k(P^\pm) \label{calU}
\end{eqnarray}
is the approximate  eigenvector from \eqref{anzatz-U}; the functions 
${V_k^\pm}$ and $u'_k$ are determined in Sections \ref{sec5.2} and \ref{sec5.3},
respectively, and the cut-off-functions and $\tau_\pm^h$ are defined in Section \ref{sec2.4}.
\ENT

According to \eqref{L'}, \eqref{a-eta}, we have
$A_k = -  u_k(P^+) / 2$ ,  $ B_k =  u_k(P^-)/2$, 
and thus $ u_k(P^-) \not= 0$ is a sufficient criterion
that the corresponding, $k$th spectral band is an interval and not
an eigenvalue of infinite multiplicity. 
The proof of Theorem \ref{th6.2} will be given in Sections \ref{sec6.2}--\ref{sec6.3},
and it needs another proof for the existence  of the gaps, 
which will be postponed to
the appendix. The proof in the appendix does not use the machinery of Sections
\ref{sec5} and \ref{sec6};  it is  an adaptation of the methods of \cite{naruta2}. 

In the following we shall also use the notation
\begin{eqnarray}
\label{use lemma}
\mu_k = \big(1+ \lambda_k+h\Lambda_k'(\eta) \big)^{-1}  \, , \ \ 
\cU_{h,k}^{\natural} = \|\cU_{h,k}; \cH^{h,\eta} \|^{-1} \cU_{h,k}  \label{6.12x}
\end{eqnarray}

\subsection{Lemma on near eigenvalues and eigenvectors.}
\label{sec6.2}
We shall need in Section \ref{sec6.3} the following operator theoretic result in the form given in \cite{BaNa}; see also \cite{BiSo} or \cite{ViLu} for more simple formulations corresponding to
$n=1$, 
$\gamma=0$ and $t = \tau$.

\BEL
\label{lem6.1}
Let ${\cB}$ be a selfadjoint, positive, and compact operator in Hilbert
space ${\cH}$ with the inner product $(\cdot,\cdot)_{\cH}$. 
If there are numbers $\mu >0$, $n \in \bbN$, and $\gamma \in (0,1/n)$,
as well as elements $\cU_1 ,\ldots , \cU_n \in {\cH}$ such that 
$\big| ( \cU_i, \cU_j)_{\cH} - \delta_{i,j} \big| \leq \gamma$ and $\|{\cB \cU_j}- \mu{\cU_j}; {\cH}
\| \leq t $ for some $t \in(0, \mu )$, then the interval 
$[\mu -\tau, \mu +\tau]$ contains at least $n$  eigenvalues of ${\cB}$,
with multiplicities counted, where $ \tau = t n^{1/2} (1 - n \gamma)^{-1/2} $.
\ENL

In this section we moreover prove the following $\eta$-independent
lower bound for the $L^2$-norm of the function $\cU_{h,k}$.

\BEL
\label{lem6.2} 
There exists a constant $C = C (k, \varpi, \theta) > 0$ such that
\begin{equation}                                                      
\label{lem1}
\Bigl|\|\cU_{h,k}; L^2 (\varpi_h) \|^2-1\Bigr|\leq 
C h \ \ \ \forall \eta \in [0,2 \pi]\,.
\end{equation}
\ENL

Proof.
First of all note that by \eqref{ortho-u}, 
$
\int_{\varpi(h)}   | u_k \circ a_h |^2 dx= a_h^{-2} ,
$
hence,  
\bea
& &  \big| \|\cX^h {u_k \circ a_h};  L^2 (\varpi_h) \|^2-1 \big|
\leq \int\limits_{\varpi_h \setminus \varpi(h)}  \!\!\!\! \big( | \cX^h u_k 
\circ a_h |^2 -|    u_k \circ a_h |^2 \big) dx  \leq Ch^3 
\label{RST1}
\eea
where the last estimate uses the facts that  $u_k$
is a bounded function near $P^\pm$, by Lemma \ref{lem2.5},
and the volume of the set $\varpi_h \setminus \varpi(h) \ni P^\pm$ is $O(h^3)$.

Denoting $\widetilde{\cU}_{h,k}=\cU_{h,k}- \cX^h {u_k \circ a_h}$, \eqref{RST1} and 
the Cauchy-Schwartz inequality yield 
\bea
& & \Big|\|\cU_{h,k};  L^2 (\varpi_h) \|^2-1)\Big| 
\leq \big|2(\cX^h {u_k \circ a_h},\widetilde{\cU}_{h,k})_{\varpi_h}\big|+\|\widetilde{\cU}_{h,k};
 L^2 (\varpi_h) \|^2 +C h 
\rowleq 
C' \|\widetilde{\cU}_{h,k}; L^2 (\varpi_h) \|^2 + C'h. \nonumber 
\eea 
Thus it is enough to prove that $\|\widetilde{\cU}_{h,k};  L^2 (\varpi_h) \| \leq C h^{1/2}$. 
By \eqref{calU}, $\|\widetilde{\cU}_{h,k}; L^2 (\varpi_h) \|$
can be bounded from above by the sum of the norms (cf. \eqref{4.77zz})
\begin{eqnarray}
\label{52}
&&\|\mathcal{X}^h_\pm \, 
u_k(P^\pm);  L^2 (\varpi_h) \|, \\
\label{53}
&&\|\chi_\pm \big(  {V_k^\mp}\circ \tau_\pm^h
-u_k(P^\pm)\big); L^2 (\varpi_h) \|,\\
\label{54}
&&h\|\mathcal{X}^hu_k' \circ a_h ; L^2 (\varpi_h) \| .
\end{eqnarray}
We complete the proof by estimating these separately. 
The  bound $ C h^{3/2} $ for \eqref{52} follows immediately by  recalling  that the supports of 
the functions $\mathcal{X}_\pm^h$ are contained in  balls of 
radius $O(h)$: the $L^2$-integral is taken over a volume bounded by $Ch^3$.

To treat \eqref{53} we recall the relation \eqref{V} and use
\eqref{5.24afNEW}:
\bea
& & \big\|\chi_\pm W_k^\mp  \circ \tau_\mp \circ  \tau_\pm^h  ;L^2(\varpi_h) \big\|^2 
\rowleq 
\int\limits_{\varpi_h }
\frac{1}{1+ h^{-2} |x- P^\pm - hP^\mp|^2 }  dx
\leq 
h^3   \!\!\! \int\limits_{B(0, 2 dh^{-1} ) }   \!\!\!
\frac{1}{1+ |\xi|^2 }  d\xi
\leq C' h^{2} ,\label{oldthing}
\eea
where the coordinate change $h^{-1} ( x - P^\pm - hP^\pm)  \mapsto \xi$  was used
and $d > 0$ denotes the diameter of the set $\varpi_h$. In the same way, using  \eqref{5.82} and  the change $h^{-1} ( x - P^\pm)   \mapsto \xi$,
we also get $\|\chi_\pm \widetilde{W}_k^\mp  \circ   \tau_\pm^h  ;L^2(\varpi_h) \|^2 \leq Ch^3$.  These estimates and \eqref{V} imply that 
\eqref{53} is bounded by $Ch$.

Finally,  $u_k'\in L^\infty (\varpi_h)$ by Remark \ref{rem3.6}, which yields the bound
$Ch$ for  \eqref{54}. 
 \ \ $\Box$

\subsection{Proof of Theorem \ref{th6.2}.}
\label{sec6.3}
To prove Theorem \ref{th6.2} we need to  prove the main estimate   \eqref{6.12}, or, 
\eqref{1.99}. For simplicity of presentation we assume here that the eigenvalue
$\lambda_k$ is simple. Multiple eigenvalues could be treated using the above formulation
of Lemma \ref{lem6.1}, but we leave the details to the reader. 
So, let us take in Lemma \ref{lem6.1} $n=1$,  $\gamma=0$ and $t = \tau$; in 
addition, we take for $\cB$   the operator 
$\cB^{h,\eta} :  \cH^{h,\eta} \to \cH^{h,\eta}$, and  $\mu = \mu_k =
(1+\lambda_k+h \Lambda_k'(\eta))^{-1}$, and $\cU_1 = \cU_{h,k}^\natural$, see \eqref{3.10} and \eqref{6.12x}. 
We are going to  prove that 
\bea
\tau = \|\cB^{h,\eta}\cU_{h,k}^\natural- \mu_k \cU_{h,k}^\natural;\cH^{h,\eta}\| \leq C h^{3/2}. 
\label{finales}
\eea
Then the lemma
implies that 
the operator $\cB^{h,\eta}$ has an eigenvalue $ \mu_k^h(\eta) $ such that  
\bea
| \mu_k^h(\eta) -\mu_k |\leq C h^{3/2} \ , \ \mbox{or} \ \ 
 | \Lambda_{j(k)}^h(\eta) - ( \lambda_k + h \Lambda_k'(\eta) )  | \leq C_k h^{3/2}  
, \nonumber
\eea 
where $ \Lambda_{j(k)}^h :=  (\mu_k^h)^{-1} - 1 $ is some eigenvalue of the model
problem. However, since $\lambda_k$ is a simple eigenvalue of the limit problem,  the result \eqref{gappi} of the Appendix implies that the only eigenvalue of the model problem, which is near $\lambda_k$, must be $\Lambda_k^h(\eta)$, i.e., $\Lambda_{j(k)}^h(\eta) =\Lambda_k^h(\eta)$. 
Thus,  \eqref{6.12} follows.

By \eqref{3.10},
\bea
& & \|\cB^{h,\eta}\cU_{h,k}^\natural- \mu_k \cU_{h,k}^\natural;\cH^{h,\eta}\|=
\sup_{Z} \langle \cB^{h,\eta}\cU_{h,k}^\natural-
\mu_k 
\cU_{h,k}^\natural,Z \rangle_{\cH^{h,\eta}}
\roweq
\sup_{Z} \Big( 
(\cU_{h,k}^\natural, Z)_{\varpi_h} 
- \mu_k ( 
\cU_{h,k}^\natural,Z )_{\varpi_h} 
- \mu_k(AD(\nabla_x)\cU_{h,k}^\natural,D(\nabla_x) Z)_{\varpi_h} \Big)
\roweq
\|\cU_{h,k};\cH^{h,\eta}\|^{-1} \mu_k \sup_{Z} |\cA ( \cU_{h,k},Z)| \ , \nonumber 
\eea
where supremum is taken over all $Z\in \mathcal{H}^{h,\eta}$ with $\|Z;\cH^{h,\eta}\|=1$
and 
\bea
& & \cA (\cU_{h,k},Z) = (\lambda_k+h\Lambda_k'(\eta))(\cU_{h,k},Z)_{\varpi_h}
- \big( AD(\nabla_x)\cU_{h,k},D(\nabla_x) Z \big)_{\varpi_h} . 
\label{5.55x}
\eea
By Lemma \ref{lem6.2},   $\|\cU_{h,k};\cH^{h,\eta}\| \geq C $,
and moreover, $\mu_k \leq 1$. So we see that the proof of 
\eqref{finales} and thus  of Theorem \ref{th6.2}
will be completed by showing that 
\bea
\sup_Z  |\cA ( \cU_{h,k},Z)| \leq Ch^{3/2} . \label{finalest}
\eea

\BEL
\label{lem7.1}
We have 
\bea
& & \cA\big( \cX^h( u_k \circ a_h  + h u_k' \circ a_h) + 
\cX_\pm^hu_k(P^\pm)
 , Z\big) 
\roweq
- h(F_1 \circ a_h,Z)_{\varpi (h)} - h(G \circ a_h,Z)_{\partial \varpi (h)} + O(h^{3/2})
\label{7.10} 
\eea
for all $Z \in \cH^{h,\eta}$,
where  $F$ and $G$ are defined in  \eqref{prob-U'-1}-\eqref{prob-U'-2} and 
we denote 
\bea
- F_1= - F + \big( -2a \lambda_k + \Lambda_k'(\eta) \big) u_k = (L_x-\lambda_k)\sum_\pm \chi_\pm 
 {\cW_k^\mp} \circ \tau_\pm . \nonumber 
\eea 
\ENL

Proof. Due to the cut-off function $\cX^h$ we can write the left hand side of 
\eqref{7.10} as 
\bea
& & 
\widetilde \cA (f,Z) := 
(\lambda_k+h\Lambda_k'(\eta))(f,Z)_{\varpi(h)}
- \big( AD(\nabla_x)f,D(\nabla_x) Z \big)_{\varpi(h)}
\label{5.55x4}
\eea
where $f =  u_k \circ a_h  + h u_k' \circ a_h$.
We first prove that 
\bea
\Big| \widetilde \cA \big( \cX_\pm ^h( u_k \circ a_h 
-u_k(P^\pm) + h u_k' \circ a_h) , Z \big)
\Big| \leq Ch^{3/2} . \label{7.12}
\eea
Since the mean value theorem and Lemma \ref{lem2.5} 
imply $|{u_k \circ a_h} (x) - u_k(P^\pm)|
\leq Ch$ in the set supp\,$\cX^h_\pm$, we get
\bea
& &  
\big|(\mathcal{X}^h_\pm( u_k \circ a_h    - u_k(P^\pm)),Z)_{\varpi(h)}\big|
\rowleq  
\Big(  \!\!\!\!\! \int\limits_{{\rm supp} \, \mathcal{X}^h_\pm \cap \varpi(h) }  \!\!\!\!\!
|{u_k \circ a_h}-u_k(P^\pm)|^2 dx
\Big)^{1/2} 
\|Z;\cH^{h,\eta}\| 
\rowleq 
C \Big( \!\!\!\!\! \int\limits_{{\rm supp} \, \mathcal{X}^h_\pm }
 \!\!\!\!\! h^2 dx \Big)^{1/2}  
 \leq Ch^{5/2}  .\nonumber 
\eea 
Moreover, since $| D( \nabla_x) \mathcal{X}^h_\pm (x) | \leq
|\nabla_x \cX_\pm ^h (x)|\leq C/h$ for all $x$, 
\bea
& & \big|(AD(\nabla_x)\mathcal{X}^h_\pm(u_k(P^\pm)-{u_k \circ a_h})),
D(\nabla_x)Z)_{\varpi(h)}\big|
\rowleq 
c\|Z;\cH^{h,\eta}\| \bigg(
\Big(  \!\!\!\!\! \int\limits_{{\rm supp} \,  \mathcal{X}^h_\pm \cap \varpi(h)}  \!\!\!\!\!
\big|\nabla_x({u_k \circ a_h}-u_k(P^\pm))\big|^2 dx
\Big)^{1/2}
\rowpl 
\Big( \!\!\!\!\! \int\limits_{{\rm supp}  \, \nabla_x \mathcal{X}^h_\pm}  \!\!\!\!\!
\big|[D(\nabla_x),\mathcal{X}^h_\pm]({u_k \circ a_h}(x)-u_k(P^\pm))\big|^2 dx \Big)^{1/2} \bigg)
\rowleq 
C\big((h^3)^{1/2}+(h^3(h^{-2}h^2))^{1/2} \big) \leq C'h^{3/2},  \nonumber 
\eea
where in the commutator, the function $ \mathcal{X}^h_\pm $
is understood as the multiplication  operator with this function, hence,
$[D(\nabla_x),\mathcal{X}^h_\pm] = 
D(\nabla_x)\mathcal{X}^h_\pm -  \mathcal{X}^h_\pm D(\nabla_x)$
is just a multiplication by the  smooth function $D( \nabla_x) \mathcal{X}^h_\pm$ 
with support contained  in supp\,$ \nabla_x 
\mathcal{X}^h_\pm $.

To estimate the same expressions with $h\mathcal{X}^h_\pm(x)  u_k'\circ a_h $ we  
take into account  Remark \ref{rem3.6}. We get, by the boundedness of $ u_k'$,
\begin{equation*}
\left|(h\mathcal{X}^h_\pm  u_k'\circ a_h , Z)_{\varpi_h}\right|\leq 
h\Big( \!\!\!\!\! \int\limits_{{\rm supp} \, \mathcal{X}^h_\pm \cap \varpi(h)} \!\!\!\!\!
|  u_k'\circ a_h |^2 dx \Big)^{1/2}\leq C h^{5/2}.
\end{equation*}
and by the $L^3$-estimate for $ \nabla_x u_k'$,
\bea
& & \left|(AD(\nabla_x) (h\mathcal{X}^h_\pm  u_k'\circ a_h), D(\nabla_x) Z)_{\varpi_h}\right|
\rowleq ch \Big( \!\!\!\!\! \int\limits_{{\rm supp}  \,  
\mathcal{X}^h_\pm\cap \varpi(h)} \!\!\!\!\!
|\nabla_x  u_k'\circ a_h |^2 dx +
h^{-2} \!\!\!\!\! \int\limits_{{\rm supp} \,  \mathcal{X}^h_\pm \cap \varpi(h)}\!\!\!\!\!
|  u_k'\circ a_h|^2 dx \Big)^{1/2} 
\rowleq
ch \bigg( \Big(  \!\!\!\!\!
\int\limits_{{\rm supp}  \, \mathcal{X}^h_\pm\cap \varpi(h)} \!\!\!\!\!
 1^{3}dx\Big)^{1/3} 
 \Big(  \!\!\!\!\! \int\limits_{{\rm supp}  \,  \mathcal{X}^h_\pm \cap \varpi(h)}  \!\!\!\!\!
|\nabla_x u_k' \circ a_h(x)|^{3} dx \Big)^{2/3} + c'h \bigg)^{1/2}\leq  C h h^{1/2}.  \nonumber 
\eea
These arguments prove \eqref{7.12}. Using that and
also \eqref{4.77zz}, the left hand side of \eqref{7.10} equals
\bea
 \widetilde \cA \big( u_k \circ a_h  + h u_k' \circ a_h  , Z \big)
 + O(h^{3/2}). \label{5.kpl}
\eea

Next, since $u_k$ is a solution of the problem \eqref{prob-3-1}--\eqref{prob-3-2} 
we get 
\bea
& & (AD(\nabla_x){u_k \circ a_h},D(\nabla_x)Z)_{\varpi (h)}=\lambda_k a_h^2({u_k \circ a_h},Z)_{\varpi (h)} 
\nonumber
\eea
and thus also the estimate
\bea
\Big|(AD(\nabla_x){u_k \circ a_h},D(\nabla_x)Z)_{\varpi (h)} -
\lambda_k(1+2ah)({u_k \circ a_h},Z)_{\varpi(h)}  \Big|\leq C h^2.
\label{5.72z}
\eea
Combining this with the definition of $\widetilde \cA$, \eqref{5.55x4},  
we get
\bea
\widetilde \cA( u_k \circ a_h, Z) = h( \Lambda_k' (\eta) - 2 a \lambda_k) (u_k,Z)_{\varpi(h)} 
+ O(h^2) .
\nonumber 
\eea
Since $u_k'$ is a solution of the problem \eqref{prob-U'-1}-\eqref{prob-U'-2}, we can write
\bea
& & h \widetilde \cA( u_k'\circ a_h, Z)
\roweq
h\Big( \big( - AD(\nabla_x) (  u_k'\circ a_h ) ,D(\nabla_x)Z
\big)_{\varpi (h)} + 
\lambda_k ( u_k'\circ a_h,Z)_{\varpi (h)} \Big)  + O(h^2) 
\roweq
- h ( F \circ a_h,Z)_{\varpi (h)} - h( G \circ a_h,Z)_{\partial \varpi (h)} + O(h^2).
\label{5.72}
\eea
Notice that here  we can commute differentiation and composition with the function $a_h$, 
for example 
$h D(\nabla_x) ( u_k'\circ a_h ) = h ( D(\nabla_x)  u_k')\circ a_h +
O(h^2) $ and so on. Now \eqref{7.10} follows from  \eqref{5.kpl}--\eqref{5.72}.
\ \ $\Box$

In view of \eqref{calU}, \eqref{5.55x}, 
\eqref{finalest}, and Lemma \ref{lem7.1}  it remains to 
show that
\bea
& & \sup\limits_Z \Big| \cA\Big( \sum_\pm  \chi_\pm\big( {V_k^\mp} \circ \tau_\pm^h 
-u_k(P^\pm) ,Z \Big)
\rowmi h(F_1 \circ a_h,Z)_{\varpi (h)} -h(G \circ a_h,Z)_{\varpi (h)}
\big)  \Big| \leq Ch^{3/2} ,
\label{5.73b}
\eea 
We  use the Green formula   for the term (cf. the second term on the right of \eqref{5.55x})
\bea
& &
-\Big(AD(\nabla_x)\chi_\pm \big( {V_k^\mp}\circ\tau_\pm^h
-u_k(P^\pm)\big) ,D(\nabla_x)Z\Big)_{\varpi(h)}   \nonumber
\eea
and write the term
inside the moduli in \eqref{5.73b} as 
\bea
& &\sum_\pm \lambda_k\Big(\chi_\pm \big({V_k^\mp}\circ\tau_\pm^h
-u_k(P^\pm)\big) ,Z\Big)_{\varpi_h} 
\rowpl 
\sum_\pm h\Lambda_k'(\eta)\Big(\chi_\pm \big({V_k^\mp}\circ\tau_\pm^h
 -u_k(P^\pm)\big),Z\Big)_{\varpi_h}
\rowmi
\sum_\pm \Big(L_x \chi_\pm ({V_k^\mp} \circ\tau_\pm^h  -
u_k(P^\pm)\big)  ,Z\Big)_{\varpi_h}  
\rowmi 
\sum_\pm \Big({N_x^h}  \chi_\pm \big({V_k^\mp} \circ\tau_\pm^h 
-u_k(P^\pm)\big)  ,Z\Big)_{\partial \varpi_h} 
\rowpl
h \Big(\sum_\pm (L_x-\lambda_k)(  \chi_\pm {\cW_k^\mp} \circ \tau_\pm ) \circ a_h ,Z 
\Big)_{\varpi (h)}
\rowpl
h  \sum_\pm \big({N_x^h} (\chi_\pm {\cW_k^\mp}\circ \tau_\pm) \circ a_h ,Z 
\big)_{\partial \varpi (h) \setminus B_h} 
\roweq
T^{(I)} + T^{(II)} +T^{(III)} +T^{(IV)} +T^{(V)} +T^{(VI)},
\label{5.73e}
\eea
where $B_h := a_h^{-1} B(P^\pm , S/2)$, cf. \eqref{Gdef}. 
 We take into account that $h \chi_\pm \circ a_h = h \chi_\pm + O(h^2)$ 
for $T^{(V)}$ and $T^{(VI)}$, and combine
\bea
S^{(I)} &:= & T^{(I)} + T^{(III)} +T^{(V)}
\roweq   
\sum_\pm \Big(  -L_x + \lambda_k) \chi_\pm \big({V_k^\mp} \circ \tau_\pm^h
-u_k(P^\pm)-h  {\cW_k^\mp}  \circ \tau_\pm   \circ a_h \big),Z \Big)_{\varpi(h)}
\rowpl
\sum_\pm \Big(( -L_x + \lambda_k ) \chi_\pm \big({V_k^\mp} \circ \tau_\pm^h
-u_k(P^\pm)\big),Z \Big)_{\varpi_h \setminus \varpi(h)}
+ O(h^2) 
\label{6.8}
\eea
(see  the remark after \eqref{5.72} for commuting the differentiation and
$a_h$).
The function  ${N_x^h} ( \chi_\pm {\cW_k^\mp}\circ \tau_\pm) \circ a_h $
vanishes in the set $\partial \varpi (h)  \cap B_h $, by  \eqref{5.24xy},
\eqref{W-}, except at the point $(a_h)^{-1}P^\pm$. We pick up
a constant $\tilde b > 2M$ ($M$ as in \eqref{5.24af}) 
such that $\varpi_h \setminus \cup_\pm B(P^\pm ,
\tilde b h ) \subset \varpi(h)$ for all $h$. Noticing that 
$ (a_h)^{-1}P^\pm \in B(P^\pm ,\tilde b h ) \subset B_h $ for small $h$,  
we can  write
\bea
 T^{(IV)} + T^{(VI)}  & =& 
\sum_\pm  \Big({N_x^h} \chi_\pm \big({V_k^\mp} \circ \tau_\pm^h
-u_k(P^\pm)-h{\cW_k^\mp} \circ \tau_\pm  \circ a_h \big),Z \Big)_{\partial \varpi(h)
\setminus B(P^\pm ,\tilde b h ) } 
\rowpl
\sum_\pm \Big({N_x^h} \chi_\pm \big({V_k^\mp} \circ \tau_\pm^h
-u_k(P^\pm),Z \Big)_{\partial \varpi_h \cap B(P^\pm ,\tilde b h )} + O(h^2) 
 =: S^{(II)} .  \nonumber 
\eea
Thus, \eqref{5.73e} equals
$O(h^2)$ plus 
\bea
\sum_\pm \big(  S^{(I)} + T^{(II)} + S^{(II)}   \big) . \label{6.9}
\eea

$1^\circ. a)$ We treat the first inner product of $S^{(I)}$, \eqref{6.8}.
We remark that $V_\mp \circ \tau_\pm^h$ 
are at least $C^2$-smooth 
in the sets $B(P^\pm , Mh)$, and the suprema of these functions and their derivatives
up to the order 2 are bounded by constants independent of $h$. 
Hence, by the same argument as below \eqref{7.12}, 
\bea
& & \big| \lambda_k\big(\chi_\pm \big({V_k^\mp}\circ\tau_\pm^h
-u_k(P^\pm)\big) ,Z\big)_{B(P^\pm , Mh)}  \big| \leq Ch^{3/2} 
\label{7.20}
\eea
and the same estimates hold even if $\lambda_k $ is replaced by $L_x$. 
Moreover, the term $L_x \cW_k^\pm $ vanishes due to \eqref{cW} and \eqref{5.24},
and \eqref{5.24a} implies
\bea
& & h\big| \lambda_k\big(\chi_\pm {\cW_k^\mp}\circ\tau_\pm
,Z \big)_{B(P^\pm , Mh) }  \big| \leq 
h \Big(  \int\limits_{B(P^\pm , Mh) } \frac{C}{|x- P^\pm|^2 } dx 
\Big)^{1/2} \leq Ch^{3/2} \label{7.22}
\eea 
As a consequence of the estimates \eqref{7.20} and \eqref{7.22}, the first
inner product in \eqref{6.8} can only be taken over the set $\varpi(h)
\setminus \cup_\pm B(P^\pm , Mh)$.

The relations \eqref{5.24a},  \eqref{cW} imply  ($a=1$), 
\bea
& &  h {\cW_k^\mp}  \circ \tau_\pm   \circ a_h (x) = 
h {\cW_k^\mp} \big( (1-  h)^{-1}  x  - P^\pm \big) 
\roweq 
(1- h) h {\cW_k^\mp} \big(  x  - P^\pm + hP^\pm \big)
\roweq 
(1 -  h){\cW_k^\mp} \big(h ^{-1}( x  - P^\pm )  - P^\mp \big)  
= (1-  h) {\cW_k^\mp}\circ \tau_\mp \circ \tau_\pm^h (x)  . \label{6.8am}
\eea
Moreover, if $|x - P^\pm| \geq Mh$, then 
$| \tau_\mp \circ\tau_\pm^h (x)| \geq M/2$, and by \eqref{cW},
\eqref{W-},  and the choice of $M$ in \eqref{5.24af} 
 we have $\cW_k^\pm \circ \tau_\mp \circ \tau_\pm^h (x) = 
W _k^\pm \circ \tau_\mp \circ \tau_\pm^h (x)$. Hence, by \eqref{V},
\bea
& & \big\|\chi_\pm  \big({V_k^\mp} \circ \tau_\pm^h-u_k(P^\pm)-
h  {\cW_k^\mp} \circ \tau_\pm  
\circ a_h \big);L^2(\varpi(h)\setminus  \cup_\pm  B(P^\pm , Mh)) \big\|^2 
\rowleq
\big\|\chi_\pm  \big({V_k^\mp} \circ \tau_\pm^h-u_k(P^\pm)-
(1-  h) {W_k^\mp}\circ \tau_\mp \circ \tau_\pm^h  \big);L^2(\varpi(h)) \big\|^2 
\rowleq
 \big\|\chi_\pm   {\widetilde{W}_k^\mp}   \circ \tau_\pm^h
 ;L^2(\varpi(h)) \big\|^2
+ 
 h \big\|\chi_\pm  {W_k^\mp} \circ \tau_\mp  \circ \tau_\pm^h  ;L^2(\varpi(h)) \big\|^2 
\label{6.8a}
\eea
Here, the argument in \eqref{oldthing} and below it yields 
\bea
& & \big\|\chi_\pm {\widetilde{W}_k^\mp}  \circ  \tau_\pm^h  ;L^2(\varpi(h)) \big\|^2 \leq 
\big\|\chi_\pm {\widetilde{W}_k^\mp}   \circ  \tau_\pm^h  ;L^2(\varpi_h) \big\|^2 
\leq C h^{3} . \label{6.11a}
\eea
Also, \eqref{oldthing} implies that the second term
on the right of \eqref{6.8a} is bounded by $Ch^3$. 

There remain the terms with $L_x$, to be evaluated on $\varpi \setminus
B(P^\pm , Mh)$.  
We first observe that due to \eqref{3.21}, \eqref{5.37c}, 
\eqref{5.24}, \eqref{cW}, 
\bea
& & L_x \chi_\pm \big( {V_k^\mp}  \circ \tau_\pm^h 
-u_k(P^\pm)-h{\cW_k^\mp} \circ \tau_\pm   \circ a_h \big)
\roweq
\widetilde L_x^\pm
\big( {V_k^\mp}  \circ \tau_\pm^h-u_k(P^\pm)
-h{\cW_k^\mp}  \circ \tau_\pm \circ a_h\big), \label{5.95y}
\eea
where $\widetilde L_x^\pm$ 
is a first order differential operator having as  coefficients bounded, 
non-constant functions with supports contained in $B(P^\pm,R) \setminus \bar B(P^\pm,S)
$  for the constants $0 < S  < R < 1/2$, see \eqref{4.77z}; notice that 
$L_x {\cW_k^\mp}  \circ \tau_\pm \circ a_h = 0$ in spite of the scaling $a_h$, since
$L_x$ only contains second order terms.  
Moreover, again by \eqref{6.8am}--\eqref{6.8a},  
\bea
 {V_k^\mp}  \circ \tau_\pm^h-u_k(P^\pm)
-h{\cW_k^\mp}  \circ \tau_\pm \circ a_h
= h {W_k^\mp}  \circ \tau_\mp \circ \tau_\pm^h
+ {\widetilde{W}_k^\mp}  \circ \tau_\pm^h
 \nonumber
\eea
so that using  \eqref{5.24afNEW}, \eqref{5.82} including the gradient
estimates, 
the argument similar to \eqref{oldthing} again yields 
\bea
& & \Big\Vert L_x \chi_\pm \big({V_k^\mp}\circ \tau_\pm^h-u_k(P^\pm)-h{W_k^\mp} 
\circ \tau_\pm  \circ a_h \big); L^2(\varpi(h)) \Big\Vert^2
 \leq C'' h^{3} . \label{6.55}
\eea

$1^\circ. b)$
 We estimate the inner product over $\varpi_h \setminus \varpi(h)$ in  
$S^{(I)}$, \eqref{6.8}. First, the term with $-L_x$ vanishes, since $V$ 
satisfies \eqref{3.21} and the supports of the coefficients of the operator
$\widetilde L_x^\pm$ do not intersect the set  $\varpi_h \setminus \varpi(h)$,
see the remarks after \eqref{5.95y}. Second, the estimate for 
the term with $\lambda_k$ is already contained in \eqref{7.20}. 
From this, \eqref{6.8} and \eqref{6.8am}--\eqref{6.55} we obtain
\bea
|S^{(I)}| \leq C h^{3/2} .  \label{6.55r}
\eea

$2^\circ.$
The estimate for  $T^{(II)}$ 
can be done as the estimate for \eqref{53} in 
\eqref{oldthing}, and it yields that 
$T^{(II)} \leq C h^3$.


$3^\circ.$
Finally,  to estimate the term $S^{(II)}$ 
we write it as the sum over $\pm$ of the terms
\bea
& & \big(  \chi_\pm {N_x^h} ({V_k^\mp}\circ \tau_\pm^h-u_k(P^\pm)-
h{\cW_k^\mp} \circ \tau_\pm  \circ a_h ),Z \big)_{\partial \varpi (h)
\setminus B(P^\pm, \tilde b h)}
\rowpl
\big(  
\widetilde N_x^h \big({V_k^\mp}\circ \tau_\pm^h-u_k(P^\pm)-
h{\cW_k^\mp} \circ \tau_\pm  \circ a_h ),Z \big)_{\partial \varpi (h)
\setminus B(P^\pm, \tilde b h) } \big) 
\rowpl
\Big({N_x^h} \chi_\pm \big({V_k^\mp} \circ \tau_\pm^h
-u_k(P^\pm),Z \Big)_{\partial \varpi_h \cap B(P^\pm, \tilde b h)}
\label{6.58}
\eea
where $ \widetilde N_x^h$  is  a bounded, smooth multiplier  with support in  $A^\pm$.

$3^\circ.a)$
To estimate the first term in \eqref{6.58} 
we observe that  the functions ${V_k^\pm} $ satisfy  \eqref{3.21k} on the boundary 
$\partial \Omega^\pm$ with $\theta \times \{ 0 \}$ excluded.  
Moreover, by \eqref{5.24xy}, \eqref{cW}, 
$N_x^h {\cW_k^\mp}  \circ \tau_\pm  \circ a_h   $ vanishes
in the set $\partial \varpi (h) \setminus B(P^\pm, \tilde b h) 
= \partial \varpi (h) \setminus B(P^\pm, \tilde b h) $, except possibly
outside  a set $B( P^\pm , \hat b) $ for some constant $\hat b>0$
(where the function $\tau_\pm \circ a_h $ does not map $\partial \varpi(h)$ 
into  the plane $\{ \xi_3 =0 \} $). Combining these observations we deduce
that the inner product can be taken only over the set $\partial \varpi(h) \setminus  B(P^\pm, \hat b )$. In this set we have, by the argument
\eqref{6.8am} and the remark after it, 
\bea
h{\cW_k^\mp} \circ \tau_\pm  \circ a_h  = 
(1-h) \cW_k^\mp \circ \tau_\mp \circ \tau_\pm^h =
(1-h) W_k^\mp \circ \tau_\mp \circ \tau_\pm^h , \nonumber
\eea
hence, proceeding as in \eqref{6.8} we  obtain
\bea
 & &\Big| 
\!\!\!\!\!\!\!\!\!\!\!\!
\int\limits_{\partial \varpi(h) \setminus  B(P^\pm, \hat b )}
 \!\!\!\!\!\!\!\!\!\!\!\!
\chi_\pm {N_x^h} \big({V_k^\mp}\circ \tau_\pm^h-u_k(P^\pm)-
h{\cW_k^\mp} \circ \tau_\pm  \circ a_h \big) Z  ds_x
\Big|
\rowleq C \Big( 
 \!\!\!\!\!\!\!\!\!\!\!\!
\int\limits_{\partial \varpi(h)  \setminus  B(P^\pm, \hat b )}
 \!\!\!\!\!\!\!\!\!\!\!\!
 \chi_\pm \big|{N_x^h}
\big(   {\widetilde{W}_k^\mp}   \circ \tau_\pm^h 
+ 
h {W_k^\mp} \circ \tau_\mp  \circ \tau_\pm^h  \big) \big|^2  ds_x \Big)^{1/2}
\Big(  \!\!\!\!\!\!\!\!\!\!\!\!
\int\limits_{\partial \varpi(h)  \setminus  B(P^\pm, \hat b ) }
 \!\!\!\!\!\!\!\!\!\!\!\! |Z|^2  ds_x \Big)^{1/2}   .
 \label{6.63}
\eea
Here, we use  \eqref{5.24afNEW}, \eqref{5.82}, and the same 
argument as in \eqref{oldthing} to estimate  
the first factor on the right by the square root of  a constant times 
\bea
& & 
 \!\!\!\!\!\!\!\!\!\!\!\!
\int\limits_{\partial \varpi(h)  \setminus  B(P^\pm, \hat b ) }
 \!\!\!\!\!\!\!\!\!\!\!\!
\Big|( \nabla_x {\widetilde{W}_k^\mp} ) \Big(\frac{x- P^\pm }{h}\Big)  \Big|^2 ds_x
+ 
h^2 
 \!\!\!\!\!\!\!\!\!\!\!\!
\int\limits_{\partial \varpi(h)  \setminus  B(P^\pm, \hat b ) }
 \!\!\!\!\!\!\!\!\!\!\!\!
\Big| {W_k^\mp} \Big(\frac{x- P^\pm - hP^\mp}{h}\Big)  \Big|^2 ds_x
\rowleq
 \!\!\!\!\!\!\!\!\!\!\!\!
\int\limits_{\partial \varpi(h)  \setminus  B(P^\pm, \hat b ) }
 \!\!\!\!\!\!\!\!\!\!\!\!
 \frac{C }{1 + h^{-6}|x- P^\pm|^6} \, ds_x
+  h^2 \!\!\!\!\!\!\!\!\!\!\!\!
\int\limits_{\partial \varpi(h)  \setminus  B(P^\pm, \hat b ) }
 \!\!\!\!\!\!\!\!\!\!\!\!
 \frac{C }{1 + h^{-4}|x- P^\pm - hP^\mp|^4} \, ds_x \leq ch^3 .
\nonumber 
\eea
The second factor on the right hand side of \eqref{6.63} 
is bounded by $C\Vert Z ; \cH^{h,\eta} \Vert$, by the Sobolev embedding 
theorem. We get  the bound $Ch^{3/2}$ for the first term in \eqref{6.58}.

$3^\circ.b)$ As for the second term in \eqref{6.58}, 
the integrand is  again supported outside the sets $B( P^\pm , S)$.
By the same argument
as in $a)$, this term is bounded by  $Ch^{3/2}.   $

$3^\circ.c)$
To evaluate the third term in \eqref{6.8} we first observe that 
the support of the function $\widetilde  N_x^h$ does
not intersect the integration domain  $\partial \varpi_h  \cap  B(P^\pm, \tilde b h) $.
Thus, on the integration domain,
\bea
 N_x^h \chi_\pm \big( V_k^\mp \circ \tau_\pm^h -u_k(P^\pm)  \big)
 =  N_x^h( V_k^\mp \circ \tau_\pm^h). \nonumber
\eea
Moreover, the function on the right vanishes everywhere else in 
$\partial \varpi_h  \cap  B(P^\pm, \tilde b h ) $ except for the 
sets $ \theta_h \times \{ P^\pm \}$, cf.\,\eqref{3.21j}--\eqref{5.37c}
and \eqref{varpi_h}. However, due to the quasiperiodicity conditions
\bea
& & 
\sum_\pm 
\int\limits_{P^\pm \times \theta_h }
\!\!\!\!
 N_x^h \big( V_k^\mp \circ \tau_\pm^h \big) \bar Z ds_x = 0 ,
 \nonumber 
\eea
hence, the third term in \eqref{6.58} vanishes.

We have thus shown that $|S^{(II)} | \leq C h^{3/2}$.
By this, \eqref{6.55r}, and  $2^\circ$, the moduli of the expressions \eqref{6.9} 
and \eqref{5.73e} have the bound $C h^{3/2}$.
This completes the proofs of \eqref{5.73b} and thus of 
\eqref{finalest}. \ \ $\Box$

\section{Appendix: existence of spectral gaps.}
\label{secNAR1}

\subsection{Upper estimate for the bands.}
\label{secNAR1.1}

In this section we complete the paper by giving a proof,
independent of the considerations in Sections \ref{sec5}--\ref{sec6},
for the existence of spectral gaps in the essential spectrum $\sigma_{\rm ess}$
of the original problem \eqref{prob-1}, \eqref{prob-2}. More precisely,
we show that   $\sigma_{\rm ess}$, corresponding to the parameter value
$h$, has a  gap between $\Upsilon_j^h$ and $\Upsilon_{j+1}^h$, 
if $j$ is such that $\lambda_j \not= \lambda_{j+1}$ and if $h$ is small enough.  

This immediately follows from the following result:
we shall show that for all $j \in \bbN$ there exist numbers
 $ h_j > 0$ and $C_j>0$ (depending also  on $\varpi$, $A$, and $\varrho$),
such that 
\bea
\lambda_j - C_j h \leq \Lambda_j^h(\eta) \leq  \lambda_j + C_j h 
\label{gappi}
\eea
for all $h \leq h_j$ and  $\eta \in [0, 2 \pi)$. 
Since this holds for all elements $\Lambda_j^h(\eta)$ of a band
$\Upsilon_j^h$, \eqref{gappi} yields a rough estimate for the 
position of the band and thus also the desired result on gaps. 

To this end we follow and modify the argument of \cite{naruta2} to prove \eqref{gappi}.
We start with an upper estimate for 
$\Lambda_j^h(\eta) \in \Upsilon_j^h $ in terms of $\lambda_j$.

\BEL
\label{lem5.4} For all $j \in \bbN$ there exist numbers
 $ h_j > 0 $ and $C_j>0$, which depend also  on $\varpi$, $A$, and $\varrho$,
and which satisfy
\bea
\Lambda_j^h(\eta) \leq  \lambda_j + C_j h 
\label{Y5.40}
\eea
for all $h \leq h_j$ and  $\eta \in [0, 2 \pi)$. 
\ENL

Proof. 
For all $j \in \bbN$, $x \in \varpi(h)$,
we set $\cV_j^h  = (1- \cX^h) {u_j^h}  $,
where $u_j^h$ is as in \eqref{2.14},  and extend the product as $0$ for
$x \in \varpi_h \setminus \varpi(h)$ (see \eqref{varpi(j,h)}, \eqref{varpi_h},
\eqref{2.14}, \eqref{4.77x}).
The definition guarantees that the extensions become smooth and
that  the support of $\cV_j^h - {u_j^h}$  is contained in a set
$\upsilon_h$ with volume estimate
$
| \upsilon_h |  \leq C h^3 .  
$
Moreover, the bounds \eqref{2.32} hold true in $\upsilon_h$.  Hence, we
get 
\bea
& & \Vert {u_j^h} - \cV_j^h ; L^2 (\varpi_h) \Vert^2 
\leq \int\limits_{ \upsilon_h}  C_j  dx
\leq C_j' h^3 , \nonumber 
\\
& & \Vert \nabla_x {u_j^h} - \nabla_x \cV_j^h ; L^2 (\varpi_h)  
\Vert^2  
\rowleq \int\limits_{\upsilon_h }
\big( | \nabla_x {u_j^h}| + | \nabla_x \cX^h| |{u_j^h}|\big)^2  dx
\leq 
C_j( h^3  +  h  ) \leq c_j' h  ,  \nonumber 
\\  
& & \Vert {u_j^h}; L^2 (\upsilon_h) \Vert^2
\leq C_j  h^3  \ , \ \ 
\Vert \nabla_x  {u_j^h}; L^2 (\upsilon_h) \Vert^2
\leq C_j  h^3  . \label{Y5.45} 
\eea

For every $j \in \bbN$ we  now pick an $h_j$, $0 < h_j < 1 $, such that
$h_j < h_{j-1}$ and, say,
$
(1 + \tilde C_j) h_j \leq 2^{-j-3},
$
where $\tilde C_j$ is the largest of the constants $C_j$ appearing
in \eqref{Y5.45}.
Let us fix the index $j$ for the rest of the proof, and consider numbers
$h$ satisfying $h < h_j$. First, by  using the choice of the numbers $h_j$, 
classical arguments and the fact that 
the functions $u_1^h, \ldots, {u_j^h}$ form an orthonormal set in 
$L^2(\varpi(h))$ (see below \eqref{2.14}), we deduce that also the small perturbations
$\cV_{1}^h , \ldots, \cV_j^h$
are linearly independent  in $L^2(\varpi(h) )^3$ and also in  $L^2 (\varpi_h)^3$.

Let now $(b_p)_{p=1}^j$ be any sequence of numbers
normalized so that $\sum_{p=1}^j
|b_p|^2 = 1$, and let 
$\cW_j^h =  \sum_{p=1}^j b_p \cV_p^h$,
$W_j := \sum_{p=1}^j b_p u_p^h ,$
hence,
$\Vert W_j  ; L^2 (\varpi(h))  \Vert = 1  $.
We evaluate 
$a( \cW_j^h, \cW_j^h ;\varpi(h))$. 
First, by \eqref{lambda},  \eqref{ortho-u}, \eqref{2.14}, 
$
|a ( W_j , W_j ;\varpi(h))| \leq \lambda_j .  \nonumber
$
Moreover, 
\bea
a(\cW_j^h  ,\cW_j^h  ;\varpi(h)) & =    & a ( W_j , W_j ;\varpi(h))
+ \sum\limits_{p,q=1}^j b_p \bar b_q
\Big( 
a(\cV_p^h  - u_p^h ,  u_q^h  ;\varpi(h)) 
\rowpl
a( u_p^h  ,\cV_q^h  -  u_q^h  ;\varpi(h)) +
a(\cV_p^h  - u_p^h  ,\cV_q^h  -  u_q^h  ;\varpi(h)) \Big)  .
\label{Y5.62}
\eea
The expression for $a(\cdot, \cdot , \varpi(h))$, \eqref{1.31},
the Cauchy-Schwartz inequality, the volume bound for $\upsilon_h$,
and \eqref{Y5.45} imply
\bea
& & |a(\cV_p^h  - u_p^h  , u_q^h  ;\varpi (h) )| 
\leq 
C_{p,q}
\Vert \cV_p^h -u_p^h ;  H^1  (\upsilon_h) \Vert 
\Vert u_q^h ;  H^1   (\upsilon_h) \Vert 
\leq C_j h^2 .
\nonumber 
\eea
The third term on the right hand side of 
\eqref{Y5.62} has the same bound, but the fourth one only the bound $Ch$.  From these estimates and
$
\sum_{p,q=1}^j b_p \bar b_q 
\leq 2j$ 
we thus obtain
\bea
& & |a(\cW_j^h  ,\cW_j^h  ;\varpi(h))| 
\leq \lambda_j +
C_j  h 
\ . \label{Y5.54}
\eea

In the same vain one can estimate using \eqref{Y5.45}
\bea
& & 
( \cW_j^h , \cW_j^h )_{\varpi(h)} 
= 
(W_j,W_j)_{\varpi(h)} +
\sum\limits_{p,q=1}^j b_p \bar b_q
\Big( (\cV_p^h  - u_p^h ,  u_q^h )_{\varpi(h)}
\rowpl
( u_p^h  ,\cV_q^h  -  u_q^h )_{\varpi(h)} +
(\cV_p^h  - u_p^h  ,\cV_q^h  -  u_q^h )_{\varpi(h)} \Big) 
\geq
1 - C_j h  .
\label{Y5.55}
\eea

To apply these estimates we use
the max--min principle \cite[Th.10.2.2]{BiSo},
\bea
\Lambda_j^h(\eta) = \max\limits_{\cH_j} \inf\limits_{\cV \in \cH_j \setminus
\{0\} } \frac{a(\cV,\cV ;\varpi_h) }{\Vert \cV ; L^2 (\varpi_h) \Vert^2} ,
\nonumber 
\eea
where $\cH_j$ stands for any subspace in $H^1_{\rm per} (\varpi_h)$ of
codimension $j-1$. 
Since the sequence $(\cV_p)_{p=1}^j$ is 
linearly independent, we can find from any $\cH_j$
an element $\bfW_j^h $ of the form 
$
\bfW_j^h (x,\eta) = e^{-i\eta z } \cW^h_j (x) . \nonumber
$
Moreover, by the definition of the form $a$, \eqref{1.31}, \eqref{variat-2}, and the 
definitions above,  we have
\bea
a(\bfW_j^h,\bfW_j^h  ;\varpi_h) = 
a(\cW_j^h,\cW_j^h ;{\varpi(h)})  , \ 
\Vert \bfW_j^h , L^2 ( \varpi_h) \Vert = 
\Vert \cW_j^h, L^2 ( {\varpi(h)}) \Vert  . \nonumber
\eea
These, together with  \eqref{Y5.54}, \eqref{Y5.55},
imply
\bea
\Lambda_j^h (\eta) & \leq  & 
\frac{a(\bfW_j^h,\bfW_j^h ;\varpi_h) }{\Vert 
\bfW_j^h ; L^2 (\varpi_h) \Vert^2}
= 
\frac{a(\cW_j^h,\cW_j^h ;{\varpi(h)}) }{\Vert 
\cW_j^h ; L^2 ({\varpi(h)}) \Vert^2}
\leq 
\frac{\lambda_j + C_j  h }{1 - C_j h }
\leq
\lambda_j + C' h  . \nonumber 
\ \Box
\eea

\subsection{Lower estimate for the bands.}
\label{secNAR1.2}

We  finally prove the lower estimate for
the numbers $\Lambda_j^h(\eta)$, see \eqref{gappi}. We use the 
following Korn inequality, which was proven in \cite{naruta2}:
for all $f \in   \cH^{h,\eta}$, there holds
\bea
\Vert f; L^2 (\varpi_h \setminus {\varpi(h)}) \Vert^2 \leq
C h^2 \big( a(f,f;\varpi_h) + \Vert f; L^2 (\varpi_h) \Vert^2 \big) .
\label{Y5.22}
\eea

\BEL
\label{lem5.2}
For every $j$ there exist a constant $C_j= C_j (\varpi, A, \varrho) >0$ and a number 
$\tilde h_j >0$  such that 
$
\Lambda_j^h(\eta) \geq \lambda_j - C_j  h 
$ for all $0 < h < \tilde h_j$. 
\ENL

Proof. 
We consider the eigenvectors $U_j^h  \in \cH^{h,\eta}$ 
see \eqref{ortho-U}: we have $a(U_j^h, U_j^h; \varpi_h) = \Lambda_j^h(\eta)$.
By \eqref{Y5.40}, if $h \leq h_j$, every $\Lambda_j^h(\eta)$
can be bounded by a positive number depending only on $j$, so \eqref{Y5.22}
implies 
\bea
& & \Vert (U_j^h ; L^2 (\varpi(h)) \Vert \geq 
\Vert (U_j^h ; L^2 (\varpi_h ) \Vert- C_jh \ , \nonumber \\ 
& & |(U_j^h , U_k^h)_{\varpi(h)} - (U_j^h , U_k^h)_{\varpi_h}| \leq C_j h \label{Y5.24x}
\eea
Thus, for  small enough  $\tilde h_j > 0 $ and $h \leq \tilde h_j$, the sequence 
$U_{1}^h, \ldots , U_j^h$ remains linearly independent in $L^2 ({\varpi(h)})^3$.

We fix $j$ and $\eta$, and assume $h \leq \tilde h_j$.
The sequence $(e^{i \eta z} U_p^h )_{p=1}^j$ is still linearly 
independent in $L^2({\varpi(h)})^3$, hence, any 
subspace $\cH_j \subset L^2({\varpi(h)})^3$ of codimension $j-1$ contains
a linear combination $\sum_{p=1}^j b_p e^{i \eta z} U_p^h$ such that
$\sum_{p=1}^j |b_p|^2 = 1 .$ We denote $ \cY_j^h = \sum\limits_{p=1}^j b_p U_p^h ;$
we have $\Vert \cY_j^h ; L^2 (\varpi_h) \Vert = 1 $.

The eigenvector property of $U_p^h$ and \eqref{Lambda} imply
$ a(\cY_j^h ,\cY_j^h ;
\varpi_h )   \leq   \Lambda_j^h(\eta) . $
On the other hand
$
|a ( e^{i \eta z} \cY_j^h , e^{i \eta z} \cY_j^h ;
\varpi_h )   | = |a(\cY_j^h ,\cY_j^h ;
\varpi_h )   | , 
$
and the estimates \eqref{Y5.24x} also yield
\bea
& & \Vert  e^{i \eta z} \cY_j^h ; L^2 ({\varpi(h)}) \Vert^2   = 
 \Vert   \cY_j^h ; L^2 ({\varpi(h)}) \Vert^2   
\geq  1 - C_j h  .
\nonumber 
\eea
Hence, we obtain 
\bea
\lambda_j & \leq &\frac{a( e^{i \eta z} \cY_j^h, e^{i \eta z} 
\cY_j^h ;{\varpi(h)}) }{\Vert  
 e^{i \eta z} \cY_j^h ; L^2 ({\varpi(h)}) \Vert^2}
\leq  
\frac{a ( e^{i \eta z} \cY_j^h, e^{i \eta z} 
\cY_j^h ;\varpi_h) }{\Vert 
 e^{i \eta z} \cY_j^h ; L^2 ({\varpi(h)}) \Vert^2}
= 
\frac{a(\cY_j^h,\cY_j^h ;\varpi_h) }{\Vert 
\cY_j^h ; L^2 ({\varpi(h)}) \Vert^2}
\rowleq 
a_\eta(\cY_j^h,\cY_j^h ;\varpi_h) ( 1 + C_j h )
\leq \Lambda_j^h(\eta) (1 +   C_j h ) . \ \ \ \Box
\nonumber 
\eea

\end{document}